\documentclass[twocolumn,letterpaper]{IEEEAerospaceCLS}

\usepackage[utf8]{inputenc}
\usepackage{bm}
\usepackage{amsmath}
\usepackage{color,soul}
\usepackage{todonotes}
\usepackage{colortbl}
\definecolor{myblue}{rgb}{0, 0.23, 0.64}
\definecolor{WVUblue}{rgb}{0, 0.16, 0.33}
\usepackage{amsfonts}
\usepackage{url}

\usepackage[breaklinks]{hyperref}
\hypersetup{
    colorlinks=true,
    linkcolor=myblue,
    filecolor=magenta,
    citecolor=myblue,
    urlcolor=myblue
}

\usepackage{algorithm}
\usepackage{algpseudocode}

\usepackage{printlen} 

\usepackage{multirow}
\usepackage{makecell}
\usepackage{lscape}
\usepackage{array}
\usepackage{nicematrix}
\usepackage{graphicx}
\usepackage[version=4]{mhchem}
\usepackage[group-separator={,}]{siunitx}
\usepackage{longtable,tabularx}
\usepackage{mathtools}
\usepackage{float}
\usepackage{mwe}
\usepackage{enumitem}
\usepackage{multirow}
\usepackage{booktabs}
\usepackage{mdframed}
\usepackage{adjustbox}
\usepackage{dsfont}
\usepackage{xspace}
\usepackage{tcolorbox}
\usepackage[flushleft]{threeparttable}
\mdfsetup{
middlelinecolor=red,
middlelinewidth=2pt,
backgroundcolor=red!5}

\begin{document}
\title{Enhancing Orbital Debris Remediation with Reconfigurable Space-Based Laser Constellations}

\author{%
David O. Williams Rogers\\ 
Department of Mechanical, Materials and Aerospace Engineering\\
West Virginia University\\
Morgantown, WV, 26506, USA\\
dw00053@mail.wvu.edu
\and 
Hang Woon Lee\\ 
Department of Mechanical, Materials and Aerospace Engineering\\
West Virginia University\\
Morgantown, WV, 26506, USA\\
hangwoon.lee@mail.wvu.edu
}

\maketitle

\thispagestyle{plain}
\pagestyle{plain}

\maketitle

\thispagestyle{plain}
\pagestyle{plain}

\begin{abstract}
Orbital debris poses an escalating threat to space missions and the long-term sustainability of Earth's orbital environment. The literature proposes various approaches for orbital debris remediation, including the use of multiple space-based lasers that collaboratively engage debris targets. While the proof of concept for this laser-based approach has been demonstrated, critical questions remain about its scalability and responsiveness as the debris population continues to expand rapidly. This paper introduces constellation reconfiguration as a system-level strategy to address these limitations. Through coordinated orbital maneuvers, laser-equipped satellites can dynamically adapt their positions to respond to evolving debris distributions and time-critical events. We formalize this concept as the Reconfigurable Laser-to-Debris Engagement Scheduling Problem (R-L2D-ESP), an optimization framework that determines the optimal sequence of constellation reconfigurations and laser engagements to maximize debris remediation capacity, which quantifies the constellation's ability to nudge, deorbit, or perform just-in-time collision avoidance maneuvers on debris objects. To manage the complexity of this combinatorial optimization problem, we employ a receding horizon approach. Our experiments reveal that reconfigurable constellations significantly outperform static ones, achieving greater debris remediation capacity and successfully deorbiting substantially more debris objects. Additionally, our sensitivity analyses identify the key parameters that influence remediation performance the most, providing essential insights for future system design. These findings demonstrate that constellation reconfiguration represents a promising advancement for laser-based debris removal systems, offering the adaptability and scalability necessary to enhance this particular approach to orbital debris remediation.
\end{abstract} 

\vspace{-3em}
\tableofcontents

\section{Introduction}

Access to space serves a critical function in present society by fostering scientific discovery \cite{danzmann_asr_2003}, international cooperation \cite{nasa_iss}, and economic development \cite{whealan_SP_2019}. In addition, new space operation concepts, such as in-orbit servicing, assembly, and manufacturing, are key drivers of the expected growth in space launches and the capitalization of the space economy \cite{faa_launch}. However, present and future missions are threatened by the rapidly expanding population of orbital debris. In particular, orbital debris forces spacecraft operators to pay for collision warnings, maneuver to avoid incoming debris, or, in the gravest outcome, mission-ending collisions (MEC) \cite{colvin_nasa_2024}.

Lethal debris, that is, objects that can lead to MEC, is grouped according to its characteristic length. Large debris denotes objects with a characteristic length greater than \SI{10}{cm}. By virtue of their size, they can be tracked from the ground, allowing operators to avoid MEC for maneuverable spacecraft. Furthermore, they are identified to be the main source of new debris \cite{mcknight_AA_2021}. Conversely, small debris refers to objects with characteristic lengths between 1 to \SI{10}{cm}, identified as the main threat to operating spacecraft.

Driven by the threats posed by debris, the literature leverages space-based lasers to remediate debris. Although not the only method proposed, interest in pulsed laser ablation (PLA) systems has grown in light of two critical aspects. First, the economic benefits of remediating debris using space-based laser platforms significantly exceed their operational cost, which is identified as the lowest compared to other methods \cite{colvin_nasa_2023}.  Second, the change in debris orbit achieved from the momentum induced by the laser platform can be used to remediate small debris only \cite{fang_aa_2019,pieters_asr_2023}, or small and large using the same platform \cite{phipps_AA_2014_ladroit,phipps_aa_2016}. The change in debris orbit, and therefore, the effectiveness of the remediation mission, is influenced by the \textit{laser-to-debris (L2D) engagement scheduling}, that is, the timing of the engagement and the selection of the engaged debris.

The adoption of an optimization framework to tackle the L2D engagement schedule directly increases the remediation performance. In Ref.~\cite{williams_asr_2025}, the authors propose an integer linear programming (ILP)-based scheduler to maximize the \textit{debris remediation capacity}. This metric quantifies the constellation's ability to nudge debris, deorbit small debris, and perform just-in-time collision avoidance (JCA). Furthermore, Ref.~\cite{baker_asc_2025} compares the performance of two debris remediation missions, one adopting a random L2D engagement schedule, and one with an optimal policy obtained using reinforcement learning. Their results show that adopting an optimized L2D policy maximizes the reduction in debris periapsis radius.

The literature sheds light on the value of L2D scheduling optimization frameworks; however, current methodologies present several limitations. State-of-the-art L2D scheduling methods assume static constellation configurations, that is, platforms do not actively change their orbits during the mission planning horizon. In effect, constraining the constellation to maintain a unique configuration leads to a lack of flexibility and responsiveness against time-sensitive events such as on-orbit breakups and collisions. Even though a few time-sensitive events can be accounted for at the time of constellation configuration design, unexpected events will occur with no prospect of a rapid and effective remediation action. In addition, for a given constellation with a fixed number of platforms, the number of feasible L2D engagements will be limited as a consequence of the dynamism and spatial distribution of the debris population, owing to their relative geometries and the inability to change them. Therefore, scalable debris remediation missions require a new paradigm to enable flexible and responsive mission planning.

The concept of \textit{constellation reconfiguration}, which involves changing the constellation configuration through the platform's orbital maneuverability, allows an increase in the debris remediation capacity. In essence, enabling the configuration to dynamically change leads to flexible and responsive mission planning. For instance, time-sensitive events such as on-orbit breakups or predicted conjunction events, originally outside of the constellation's operational L2D range, can be tackled by reconfiguring to a constellation configuration capable of engaging them and, therefore, effectively mitigating the threat. Following the same rationale, reconfiguring the constellation configuration over time offers the possibility to perform ablation on debris originally outside of the initial constellation's operational range. In addition, adapting the constellation configuration to the highly dynamic debris population facilitates obtaining better L2D relative geometries, which significantly influence the magnitude and direction of the imparted $\Delta v$. Despite the benefits offered by constellation reconfiguration, the literature lacks an optimization framework for its implementation on space-based debris remediation missions. In contrast, several mathematical optimization approaches have been used to quantify their benefits in other application domains.

Complex and dynamic space missions' objectives are tackled in the literature by leveraging constellation reconfiguration optimization. For instance, a study in Ref.~\cite{pearl_jsr_2025} demonstrates that, for tracking fast-paced dynamical events, constellation reconfiguration provides improved observational coverage compared to satellite agility (the ability to perform slewing maneuvers). Furthermore, Ref.~\cite{lee_jsr_2023} proposes single-stage (\textit{i.e.,} single opportunity) reconfiguration for disaster monitoring by maximizing the observation rewards for a set of static targets, given $\Delta v$ budget constraints. Reference~\cite{lee_jsr_2024} extends the previous formulation to multistage reconfiguration and demonstrates the benefits of increasing the number of stages by maximizing a Hurricane's observation rewards. In addition to multistage reconfiguration, Ref.~\cite{pearl_arxiv_2025} proposes to integrate memory storage management and data downlink to the scheduler, further proving the benefits of constellation reconfiguration into the Earth observation satellite scheduling problem. 
 
\begin{figure}[htpb]
    \centering
    \includegraphics[width=0.87\linewidth]{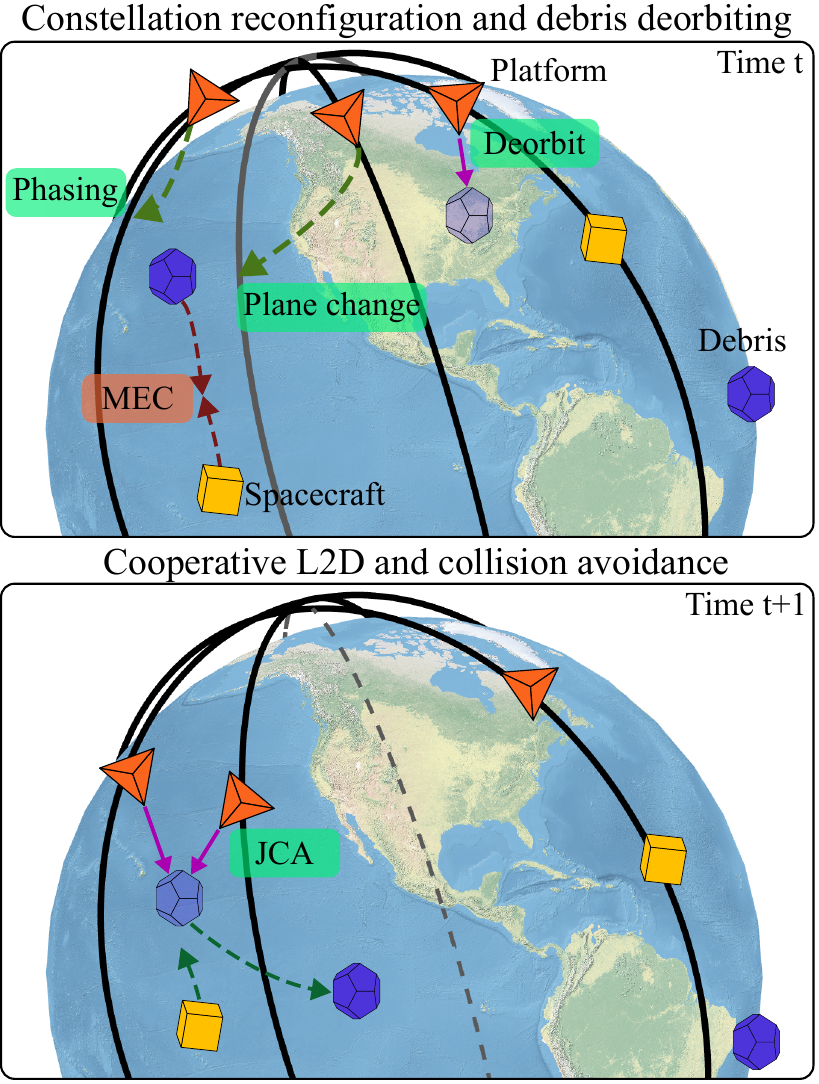}
    \caption{Proposed R-L2D-ESP's key functionalities.}
    \label{fig:overarching}
\end{figure}

Motivated by the proven advantages of constellation reconfiguration and seeking to enable flexible and responsive debris remediation using a constellation of space-based lasers, we propose:
\begin{enumerate}[leftmargin=*]
    \item The \textit{Reconfigurable L2D Engagement Scheduling Problem} (R-L2D-ESP), an optimization problem that maximizes the debris remediation capacity by determining a provably optimal schedule for constellation reconfiguration and cooperative L2D engagements.
    \item A \textit{receding horizon scheduling} (RHS) approach to address the rapid expansion in the  R-L2D-ESP solution space, resulting from its combinatorial nature. The RHS breaks the R-L2D-ESP into subproblems with shorter planning horizons and solves them recursively, leveraging information from future subproblems.
    \item A set of computational experiments that prove the value of adopting constellation reconfiguration for orbital debris remediation. In addition, sensitivity analyses are presented to show the impact of distinct mission parameters on the debris remediation capacity.
\end{enumerate}
Figure~\ref{fig:overarching} presents the key features of the R-L2D-ESP. The top block illustrates a forecasted MEC between an operational spacecraft (yellow cube) and debris (blue dodecahedron), indicated by red dashed lines, a platform performing a phasing maneuver, and another one performing a plane change maneuver, both indicated by green dashed lines, yielding a constellation reconfiguration. Further, a platform performs an L2D engagement, represented with a magenta solid line, that leads to debris being deorbited. The bottom block shows a cooperative L2D engagement to perform JCA between the spacecraft and debris.

The remainder of the paper is organized as follows. Section~\ref{sec:methodology} introduces the novel R-L2D-ESP, and Sections~\ref{sec:case_studies_1} and~\ref{sec:breakup} showcase its benefits against a non-reconfigurable constellation. Section~\ref{sec:sensitivity} presents sensitivity analyses that characterize the impact of key parameters on the debris remediation capacity. Lastly, Section~\ref{sec:conclusions} concludes the paper and proposes future research directions.

\section{Methodology} \label{sec:methodology}
The R-L2D-ESP uses a finite number of discrete sets. First, the set of time steps, or the mission planning horizon ${\mathcal{T}} \coloneqq \{0, \dots, T-1\}$, with index $t$ and cardinality $T$. Second, the set of (laser) platforms ${\mathcal{P}}\coloneqq \{0,\dots,P-1\}$, with index $p$ and cardinality $P$. Third, the set of debris ${\mathcal{D}}\coloneqq \{0,\dots,D-1\}$, with index $d$ and cardinality $D$. Each debris $d$ has, at time step $t$, a position and velocity vectors $\bm{r}_{td},\bm{v}_{td}$, respectively, that fully define the state vector $\bm{q}_{td}=[\bm{r}_{td},\bm{v}_{td}]^T$. In addition, the L2D engagements and the delivered $\Delta v$ are obtained by leveraging the mechanism of PLA. Time-expanded graphs (TEG) are adopted to model the sequential orbital transfers, during the mission planning horizon, for debris and platforms. Lastly, Table~\ref{table:fundamental_parameters} presents the parameters introduced in this section.

\begin{table}[htpb]
\renewcommand{\arraystretch}{1.3}
\caption{\textbf{Fundamental discrete sets and common parameters.}}
\label{table:fundamental_parameters}
\centering
\begin{tabular}{|c|l|}
\hline
\bfseries Symbol & \bfseries Description \\
\hline\hline
${\mathcal{T}}$  & Set of time steps (index $t$, cardinality $T$) \\
${\mathcal{P}}$  & Set of platforms (index $p$, cardinality $P$) \\
${\mathcal{D}}$  & Set of debris (index $d$, cardinality $D$) \\
$\bm{r}_{td}$ & Debris $d$ position vector at time step $t$\\
$\bm{v}_{td}$ & Debris $d$ velocity vector at time step $t$\\
$\bm{q}_{td}$ & Debris $d$ state vector at time step $t$\\
\hline
\end{tabular}
\end{table}

\subsection{Pulsed Laser Ablation}

PLA is the physical mechanism by which an object's surface material is removed using a pulsed laser. Among all the factors that dictate the PLA mechanism, the pulse energy has a critical influence \cite{vonderlinde_apss_2000}. Particularly, for nanosecond and picosecond pulse lengths, PLA induces surface material's heating, melting, and vaporization \cite{Stafe_LA_book_2014}. The laser beam's energy flux is absorbed by the vapor, heating and ionizing it, leading to plasma formation \cite{lunney_ass_1998}. Then, plasma expands laterally and perpendicularly to the surface, inducing a momentum on the target object. Among the multiple factors that characterize the momentum transfer mechanism, a key parameter is the momentum coupling coefficient $c_\text{m}$. This coefficient quantifies the amount of optical pulse energy that is transferred as momentum, defined as \cite{yuan_jap_2012}:
\begin{equation}
     c_\text{m}=\frac{I}{E}
\end{equation}
where $I$ indicates debris momentum, and $E$ the optical pulse energy.

The literature (\textit{e.g.,} Refs.~\cite{pieters_asr_2023,phipps_AA_2014_ladroit,phipps_aa_2016,soulard_AA_2014}) proposes to leverage the mechanism of PLA to induce an impulsive $\Delta v$ owed from plasma expansion, and change debris orbit. Assuming a perfect spherical debris, with uniform surface mass density $\mu$, the instantaneous $\Delta v$ delivered by a single pulse is computed as \cite{phipps_AA_2014_ladroit}:
\begin{equation}
    \Delta v = \eta_1 \frac{ c_\text{m}\varphi}{\mu}\label{eq:DeltaV}
\end{equation}
where $\eta_1$ is an efficiency coefficient that considers improper thrust direction, debris shape, and tumbling, and $\varphi$ is the range-dependent laser fluence delivered on debris surface, defined as \cite{phipps_AA_2014_ladroit}:
\begin{equation}
    \varphi = \eta_2\frac{E}{\pi}\left(\frac{2\phi}{M^2
    a\lambda u}\right)^2\label{eq:fluence}
\end{equation}
where $\eta_2$ is an efficiency coefficient that considers system losses derived from apodization and obscuration, $\phi$ is the laser system's primary mirror diameter, $M^2$ is the beam quality, $a$ is the diffraction constant, $\lambda$ is the wavelength, and $u$ is the L2D range.

Letting $\bm{r}_{tp}$ be the position of platform $p$ at time step $t$, the L2D range $u_{tpd}$ at time step $t$ between platform $p$ and debris $d$ is defined as $\|\bm{u}_{tpd}\|_2 =\|\bm{r}_{td}-\bm{r}_{tp}\|_2$. Further, the unit vector that defines the direction of the L2D is given as $\hat{\bm{u}}_{tpd}=\bm{u}_{tpd}/u_{tpd}$. By leveraging it and assuming that debris is a perfect sphere, we can now define the imparted $\Delta \bm{v}$ as:
\begin{equation}
    \Delta \bm{v} = \Delta v \hat{\bm{u}}= \eta_1 \frac{ c_\text{m}\varphi}{\mu}\hat{\bm{u}}\label{eq:DeltaV_vector}
\end{equation}

Adopting the fundamental sets outlined in Table~\ref{table:fundamental_parameters}, one can define the total $\Delta \bm{v}_{tpd}$ imparted by platform $p$ on debris $d$ at time step $t$, as a function of the number of pulses $n_{\text{L2D}}$ delivered on a single time step, as:
\begin{equation}
    \Delta \bm{v}_{tpd} = n_{\text{L2D}} \eta_1\frac{ c_{\text{m},{tpd}}\varphi_{tpd}}{\mu_{d}}\hat{\bm{u}}_{tpd}\label{eq:DeltaV_windex}
\end{equation}
The number of pulses, $n_{\text{L2D}}$, is computed as the product of the laser's pulse repetition frequency and the length of the L2D engagement. The laser is assumed to operate in a shooting-cooling mode scheme, where the shooting mode consists of the laser performing an L2D engagement, and the cooling mode denotes the period of time that the accumulated heat is removed from the laser. 

To characterize simultaneous and cooperative L2D engagements from multiple platforms on the same debris object, we adopt the $\Delta v$ vector analysis (DVA) framework proposed in Ref.~\cite{williams_asr_2025}. Specifically, the framework quantifies the cooperative and simultaneous L2D engagements at time step $t$ on debris $d$ from a subset of platforms ${\mathcal{P}}_{td}\subseteq{\mathcal{P}}$, through the summation of its vector components as defined in Eq.~\eqref{eq:DeltaV_windex}. Further, assuming that the imparted $\Delta v$ is instantaneous, debris $d$ position at time step $t$ before an L2D engagement $\bm{r}^-_{td}$, and after $\bm{r}^+_{td}$, are equal, and hence it is denoted as $\bm{r}_{td}$. Then, given its state vector $\bm{q}_{td}$ before the engagement, and the set of platforms ${\mathcal{P}}_{td}$ that perform L2D ablation on it, the post-L2D engagement state vector $\bm{q}^+_{td}$ is defined as:
\begin{equation}
    \bm{q}^+_{td}=\biggl[\bm{r}_{td},\bm{v}_{td} + \sum_{p \in {\mathcal{P}}_{td}}\Delta \bm{v}_{tpd}\biggr]^T \label{eq:state_vec_dva}
\end{equation}
Lastly, debris $d$ post-L2D engagement periapsis radius $r^+_{\text{[peri]},td}$ at time step $t$, is analytically derived from its post-L2D engagement state vector $\bm{q}^+_{td}$ as:
\begin{equation}
    r^+_{\text{[peri]},td} = \frac{ \|(\bm{r}_{td}\times\bm{v}^+_{td})\|_2^2}{\mu_{\text{Earth}}}\frac{1}{1+e^+}\label{eq:debris_periapsis}
\end{equation}
where $\mu_{\text{Earth}}$ is the Earth's gravitational parameter, $\bm{v}^+_{td}$ debris $d$ post-L2D engagement at time step $t$, and $e^+$ the post-engagement eccentricity. Lastly, Table~\ref{table:pla_parameters} summarizes all parameters introduced in this section.
\begin{table}[htpb]
\renewcommand{\arraystretch}{1.3}
\caption{\textbf{PLA parameters.}}
\label{table:pla_parameters}
\centering
\begin{tabular}{|c|l|}
\hline
\bfseries Symbol & \bfseries Description \\
\hline\hline
$c_{\text{m}}$  & Momentum coupling coefficient \\
$I$ & L2D impulse \\
$E$ & Pulse energy\\
$\eta$ & Efficiency coefficient\\
$\varphi$ & L2D fluence \\
$\mu$ & Debris surface density\\
$\phi$ & Mirror diameter\\
$M^2$ & Beam quality\\
$a$ & Difraction constant\\
$\lambda$ & Wavelength\\
$u$ & L2D range\\
$n_{\text{L2D}}$ & Number of pulses per time step\\
$\bm{q}^+_{td}$ & Debris $d$ state vector post-L2D at time step $t$\\
$e$ & Orbit eccentricity\\
\hline
\end{tabular}
\end{table}

\subsection{Orbital Transfers Using Time-Expanded Graphs}

In this paper, we adopt TEGs to model the sequence of platform and debris orbital transfers over time. In essence, TEGs are directed graphs where each vertex represents an orbital slot defined at time step $t$, and each directed edge constitutes the transfer between two orbital slots. Specifically, an orbital slot fully defines an orbit by encoding the state vector defined at time step $t$ for a platform or debris. In addition, each transfer (\textit{i.e.,} edge) has an associated weight that represents the cost or benefit of transferring from an orbital slot defined at the time step $t$ to an orbital slot defined at $t+1$.

Platforms change their orbital slots during the mission planning horizon using impulsive orbital maneuvers. Let ${\mathcal{S}}^p_{t}\coloneqq \{0,\dots,S^p_{t}-1 \}$ be the set of orbital slots for platform $p$ at time step $t$, with index $s$ and cardinality $S^p_{t}$. Particularly, at the epoch $t=0$, each platform is initialized at orbital slot $s=0$; therefore, ${\mathcal{S}}^p_{0}\coloneqq \{0\}$ for all $p$ in ${\mathcal{P}}$. Then, we can define for each platform $p$, a TEG ${\mathcal{G}}^p =({\mathcal{S}}^p, {\mathcal{E}}^p)$, where ${\mathcal{S}}^p\coloneqq \{ {\mathcal{S}}^p_{0}, \dots, {\mathcal{S}}^p_{t}, \dots,  {\mathcal{S}}^p_{T-1}\}$ denotes the set of vertices (\textit{i.e.,} orbital slots) for platform $p$ during the mission planning horizon and ${\mathcal{E}}^p$ the set of edges (\textit{i.e.,} transfers) between orbital slots. For each time step $t$ and platform $p$, we define the orbital transfer cost $c^p_{tsw} \in {\mathbb{R}}_{\ge 0}$, which encodes the required $\Delta v$ to transfer from orbital slot $s \in {\mathcal{S}}^p_{t}$ to $w \in {\mathcal{S}}^p_{t+1}$. The cost is obtained by solving an orbital boundary value problem given initial and final conditions defined by each orbital slot. However, if $s=w$, the platform is assumed to maintain its orbit and the transfer cost is $c^p_{tsw}=0$.

\begin{figure}[htpb]
    \centering
    \includegraphics[width=\linewidth]{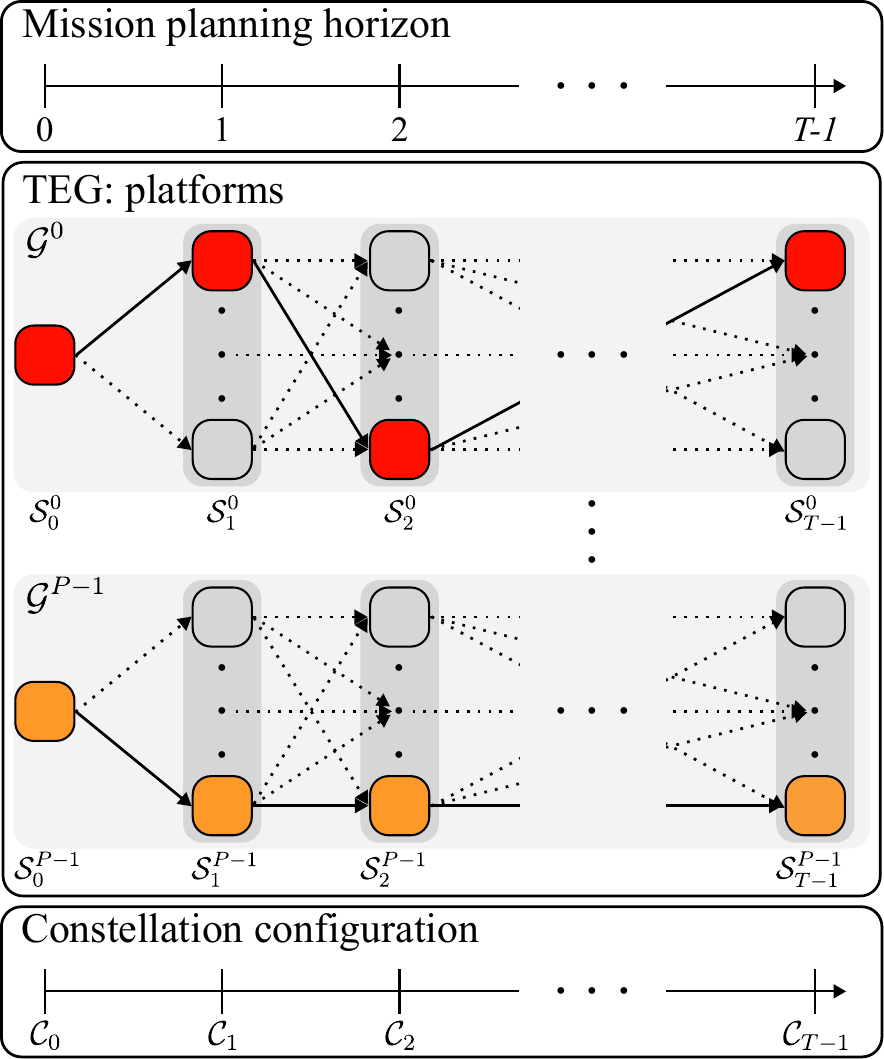}
    \caption{Illustration of TEGs for platforms $\{0,\dots,P-1\}$.}
    \label{fig:platform_teg}
\end{figure}

Figure~\ref{fig:platform_teg}, illustrates platforms $p=0$ and $p=P-1$ TEGs. At each time step $t=1$, a constellation configuration ${\mathcal{C}}_1$ is defined by the selected platform's orbital slots, represented by colored squares. Similarly, for time step $t=2$, there is an associated constellation configuration ${\mathcal{C}}_2$ defined by its corresponding orbital slots. The reconfiguration process between two constellation configurations is given by the platform's orbital transfers, indicated by solid lines.

To model candidate transfers between orbital slots and determine which ones are selected by the scheduler, we adopt binary decision variables:
\begin{equation}
    z^p_{tsw}=\begin{cases}
        1, \quad & \text{if at time step $t$, platform $p$ transfers}\\
                &\text{from orbital slot $s \in {\mathcal{S}}^p_{t}$ to $w \in {\mathcal{S}}^p_{t+1}$}\\
        0, \quad & \text{otherwise}
    \end{cases}
\end{equation}
The set of selected platforms' transfer paths determines the reconfiguration process at time step $t$ between the constellation configurations ${\mathcal{C}}_t$ and ${\mathcal{C}}_{t+1}$.

Debris transfers to a new orbit if at least one platform engages it. At time step $t$, platform $p$ is allowed to engage debris $d$ if it satisfies two conditions. The first condition, described in Eq.~\eqref{eq:cond_visibility}, enforces debris $d$ to be in the line-of-sight of platform $p$ at time step $t$. The second condition, specified in Eq.~\eqref{eq:cond_range}, enforces that the L2D range $u_{tpd}$ is greater than a lower bound $u_{\min}$, defined for operational safety, and lower than an upper bound $u_{\max}$, which encodes the distance at which the delivered fluence is lower than the ablation threshold, \textit{i.e.,} no $\Delta v$ induced.
\begin{subequations}
    \begin{alignat}{2}
        &0\le\sqrt{(r^p_{ts})^2-(R_{\oplus}+\epsilon)^2} + \notag \\
        & \qquad \qquad \qquad \sqrt{(r_{td})^2-(R_{\oplus}+\epsilon)^2} 
        -u^p_{tsd}\label{eq:cond_visibility}\\
        &u_{\min}\le u^p_{tsd} \le u_{\max}\label{eq:cond_range}        
    \end{alignat}
\end{subequations}
To encode if a platform engages debris, we propose binary decision variables:
\begin{equation}
    y^p_{tsd}=\begin{cases}
        1, \quad & \text{if at time step $t$, platform $p$ located at}\\
                &\text{orbital slot $s \in {\mathcal{S}}^p_{t}$ engages debris $d$}\\
        0, \quad & \text{otherwise}
    \end{cases}
\end{equation}

\begin{figure*}[htpb]
    \centering
    \includegraphics[width=\linewidth]{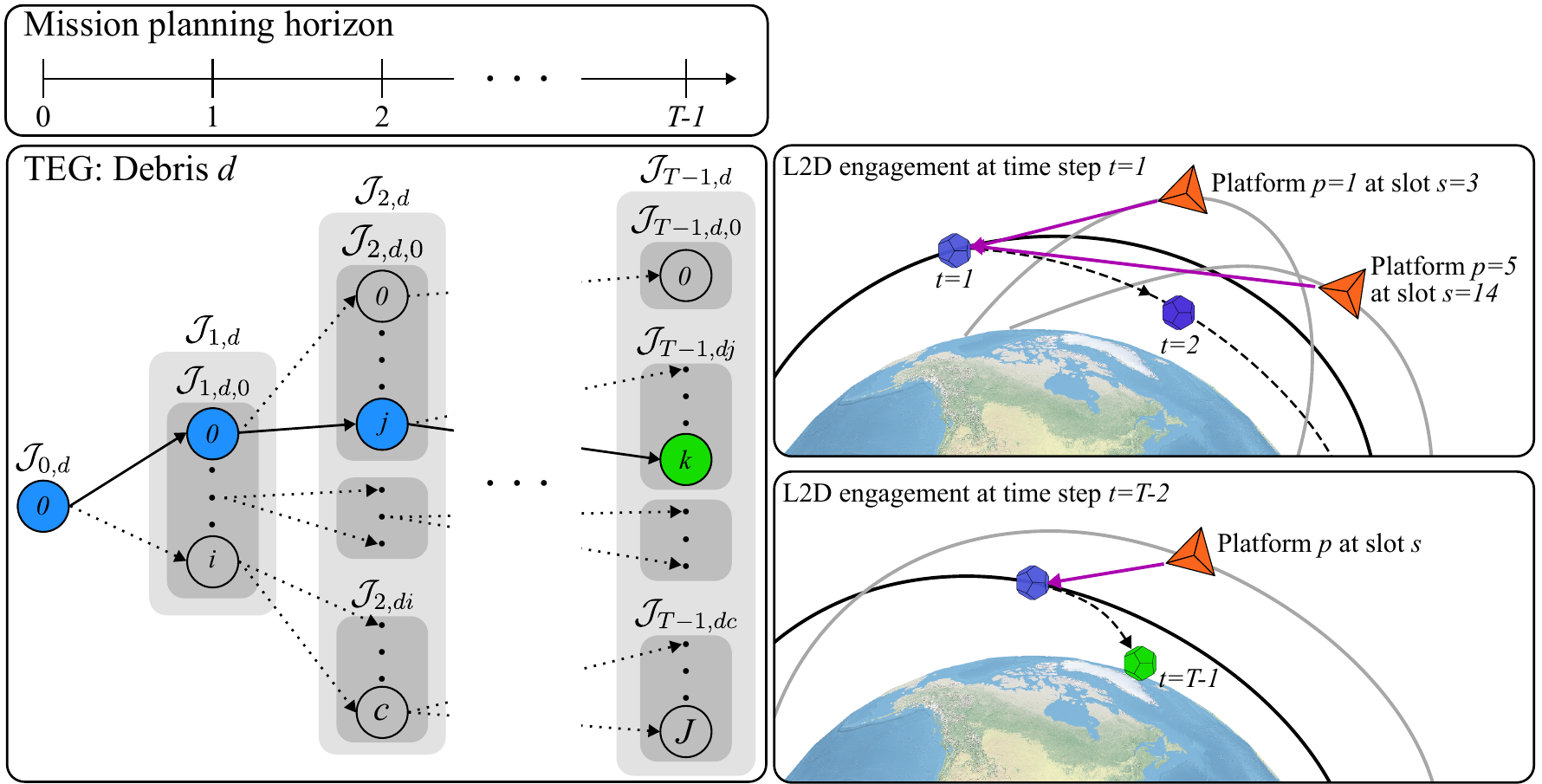}
    \caption{Illustration of debris $d$'s TEG. The circles represent orbital slots $j \in {\mathcal{J}}_{td}$, with blue indicating those selected, and green indicating that an orbital slot represents that debris is deorbited. The dashed lines represent transfer opportunities, and the solid ones those selected. The upper-right block shows a collaborative engagement at time step $t$ from the set of platforms ${\mathcal{(PS)}}_{1,dj}=\{(1,3),(5,14)\}$, encoding platforms $p=1$ and $p=5$ at orbital slots $s=3$ and $s=14$, respectively. They generate orbital slot $j \in {\mathcal{J}}_{2,d,0}$. Similarly, the right-lower block encodes at time step $T-2$ the transfer to slot $k \in {\mathcal{J}}_{T-1,d}$, which encodes deorbiting, generated by platform $p$ at orbital slot $s$ defined in set ${\mathcal{(PS)}}_{T-2,dk}=\{(p,s)\}$.}
    \label{fig:debris_teg}
\end{figure*}

Figure~\ref{fig:debris_teg} illustrates debris $d$'s TEG ${\mathcal{H}}_d=({\mathcal{J}}_d, {\mathcal{E}}_d)$. The set of TEG vertices ${\mathcal{J}}_d \coloneqq \{ {\mathcal{J}}_{0,d}, \dots, {\mathcal{J}}_{td}, \dots,  {\mathcal{J}}_{T-1,d}\}$ represents the set of debris orbital slots during the mission planning horizon. Each orbital slot $j \in {\mathcal{J}}_{td}$ has an associated set ${\mathcal{(PS)}}_{t-1,dj}\coloneqq\{(p,s),\dots,(p_{n}, s_{n})\}$ with indices $p,s$ and cardinality $n$. The set ${\mathcal{(PS)}}_{t-1,dj}$ encodes the platforms $p$ with their respective orbital slots $s$ required to engage debris $d$ at time step $t-1$ and generate orbital slot $j \in {\mathcal{J}}_{td}$. Furthermore, each new debris orbital slot has an associated parent orbital slot, owing to the tree structure of debris TEG. To this end, we define for all time steps $t$ greater than zero and all debris $d$, the set ${\mathcal{J}}_{tdi}$, with index $j$, and cardinality $J_{tdi}$, which indicates all child orbital slots $j \in {\mathcal{J}}_{tdi}$ corresponding to parent orbital slot $i \in {\mathcal{J}}_{t-1,d}$. 

The set of debris $d$'s TEG's edges ${\mathcal{E}}_d$ represents the transfer between two orbital slots. At time step $t$, debris $d$'s transfer from orbital slot $i \in {\mathcal{J}}_{td}$ to $j \in {\mathcal{J}}_{t+1,di}$ has an associated reward $R_{tdij}$, defined as:
\begin{equation}
    R_{tdij}=\Gamma_{tdij}+\gamma_{\text{[peri]},tdij}
\end{equation}
where $\Gamma_{tdij}$ is a penalty term that takes a large negative value if debris $d$'s transfer at time step $t$ from orbital slot $i \in {\mathcal{J}}_{td}$ to $j \in {\mathcal{J}}_{t+1,di}$ creates a conjunction event with any active spacecraft considered in the optimization. Particularly, the user can define a set of active spacecraft ${\mathcal{A}}$ and propagate their orbits during the mission planning horizon. If any debris post-L2D engagement orbital slot $j$ invades the conjunction ellipsoid of an active spacecraft, then $\Gamma_{tdij}=-\alpha$, where $\alpha$ is a user-defined large penalty value; therefore, discouraging the constellation from engaging that debris. Otherwise, $\Gamma_{tdij}=0$. The second term $\gamma_{\text{[peri]},tdij}$ quantifies the L2D engagement as a function of the new orbital slot's periapsis radius. Formally, $\gamma_{\text{[peri]},tdij}$ is defined as:
\begin{equation}
    \gamma_{\text{[peri]},tdij}=\begin{cases}
        (r_{\text{deorbit}}/r^+_{\text{[peri]},td})^3, \quad & \text{if periapsis radius is}\\
                    &\text{lowered} \\
        100, \quad & \text{if debris $d$ is deorbited}\\
        0, \quad & \text{if $i$ and $j$ encode the}\\
                    &\text{same orbit}\\
    \end{cases} \label{eq:reward}
\end{equation}
where each reward case represents: (1) a transfer yielding a periapsis radius closer to the deorbit radius threshold $r_{\text{deorbit}}$, (2) deorbit, and (3) no orbital transfer. If the post-engagement periapsis radius $r^+_{\text{[peri]},td}$ is lower or equal to the deorbit radius threshold $r_{\text{deorbit}}$, reward case (1) is greater than one; however, to avoid numerical issues we propose reward case (2) that assigns the same reward of 100 to a debris transfer that results in deorbit. Alternatively, the post-engagement periapsis radius can be compared with the pre-engagement periapsis radius $r_{\text{[peri]},td}$; however, this approach does not provide a sense of how close the debris object is to the deorbit condition. Further, if the L2D engagements from at least one platform increase the periapsis radius, we do not generate that orbital slot. Hence, there is no $\gamma_{\text{[peri]},tdij}$ value associated with it. 

We propose Algorithm~\ref{alg:debris_teg} to generate debris $d$'s TEG and its necessary parameters: transfer reward $\bm{R}$, the set of orbital slots ${\mathcal{J}}_d$ and ${\mathcal{J}}_{tdi}$, and the set of platforms ${\mathcal{(PS)}}_{tdj}$. The algorithm sequentially iterates over each time step $t$ and starts by adding an orbital slot $j$, which encodes that debris orbit remains unchanged; therefore, the associated reward is zero. Then, it iterates over the set of orbital slots $i \in {\mathcal{J}}_{td}$, and determines if that orbital slot encodes a deorbit condition, indicated by a state vector $\bm{q}_{td}=\bm{0}$. Otherwise, the function \texttt{GetPeriapsis} determines the periapsis radius of debris $d$ at time step $t$ using Eq.~\eqref{eq:debris_periapsis}, and the function \texttt{isL2DFeasible} generates a list of all combinations of platforms, and their orbital slots that satisfy the L2D engagement conditions defined in Eqs.~\eqref{eq:cond_visibility} and~\eqref{eq:cond_range}.

Iterating over the list obtained with \texttt{isL2DFeasible}, and using Eq.~\eqref{eq:state_vec_dva} to obtain, for time step $t$, debris $d$ post-L2D state vector $\bm{q}^+_{td}$, we compute if the new orbit has a lower periapsis radius $r^+_{\text{[peri]},td}$. If true, we store the new orbital slot, save the combination of platforms and their orbital slots that generate it, and call \texttt{GetReward} that uses Eq.~\eqref{eq:reward} to compute the transfer reward. When the reward is lower than one but greater than zero, it indicates that debris transfers to a new orbit, which is obtained by propagating the orbit with the new state vector using \texttt{PropagateOrbit}. Otherwise, the engagement leads to deorbit; therefore, the state vector for the new orbital slot is set to zero for the remainder of the mission planning horizon $\{t+1,\dots,T-1\}$.
\begin{algorithm}[htpb]
\scriptsize
    \caption{Debris $d$'s TEG generation and necessary parameters.}
    \label{alg:debris_teg}
    \begin{algorithmic}[1]
    \renewcommand{\algorithmicrequire}{\textbf{Input:}}
    \renewcommand{\algorithmicensure}{\textbf{Output:}}
        \Require ${\mathcal{T}}, {\mathcal{P}}, {\mathcal{S}}^p,{\mathcal{D}}$
        \Ensure ${\mathcal{J}}_d, R_{tdij}, {\mathcal{J}}_{tdi}, {\mathcal{(PS)}}_{tdj}$
        \For{$t \in {\mathcal{T}}\setminus\{T-1\}$}
            \State $j \gets 0$
            \For{$i \in {\mathcal{J}}_{td}$}
                \State ${\mathcal{J}}_{t+1,d} \gets {\mathcal{J}}_{t+1,d} \cup \{j\}$
                \State ${\mathcal{J}}_{t+1,di} \gets {\mathcal{J}}_{t+1,di} \cup \{j\}$
                \State $\bm{q}_{t+1,dj} \gets \bm{q}_{t+1,di}$
                \State $R_{tdij}=0$
                \If{$\bm{q}_{tdi}== \bm{0}$}
                    \State $\bm{q}_{t+1,dj} \gets \bm{0}$
                    \State $j \gets j + 1$
                \Else
                    \State $r_{\text{[peri]},td} \gets \texttt{GetPeriapsis}(\bm{q}_{tdi})$
                    \State ${\mathcal{\bar{P}}}_{tdi} \gets \texttt{isL2DFeasible}(\bm{q}_{tdi}, {\mathcal{S}}^p_t)$
                    \For{$(\bm{p}, \bm{s}) \in {\mathcal{\bar{P}}}_{td}$}
                        \State $r^+_{\text{[peri]},td} \gets \texttt{GetPeriapsis}(\bm{q}^+_{tdi})$ 
                        \If{$r^+_{\text{[peri]},td} < r_{\text{[peri]},td}$}
                            \State $j \gets j + 1$
                            \State ${\mathcal{J}}_{t+1,d} \gets {\mathcal{J}}_{t+1,d} \cup \{j\}$
                            \State ${\mathcal{J}}_{t+1,di} \gets {\mathcal{J}}_{t+1,di} \cup \{j\}$
                            \State ${\mathcal{(PS)}}_{tdj} \gets \{(p,s),\dots,(p_{n}, s_{n})\}$
                            \State $R_{tdij} \gets \texttt{GetReward}(r_{\text{[peri]},td}, r^+_{\text{[peri]},td})$
                            \If{$R_{tdij}<1$}
                                \State $\bm{q}_{t+1,dj} \gets \texttt{PropagateOrbit}(\bm{q}^+_{tdi})$
                            \Else
                                \State $\bm{q}_{t+1,dj} \gets \bm{0}$                               
                            \EndIf
                        \EndIf
                    \EndFor
                    \State $j \gets j + 1$
                \EndIf
            \EndFor
        \EndFor
    \end{algorithmic}
\end{algorithm}

Lastly, the transfer path between debris orbital slots during the mission planning horizon is encoded in binary decision variables:
\begin{equation}
    x_{tdij}=\begin{cases}
        1, \quad & \text{if at time step $t$, debris $d$'s transfers from}\\
                &\text{orbital slot $i \in {\mathcal{J}}_{td}$ to $j \in {\mathcal{J}}_{t+1,di}$}\\
        0, \quad & \text{otherwise}
    \end{cases}
\end{equation}

\subsection{The Reconfigurable Laser-to-Debris Engagement Scheduling Problem}

The R-L2D-ESP schedules the sequence of constellation reconfigurations and L2D engagements that maximize the debris remediation capacity. The problem is modeled as an ILP given the binary domain of the decision variables and the rapidly expanding TEG structures of debris and platform transfers. The key benefit of adopting ILP is that it enables us to obtain provably optimal solutions using commercial solvers such as the Gurobi Optimizer or CPLEX. We refer the reader to Appendix~\ref{app:nomenclature} for a review of the nomenclature used in this section.

To maximize the constellation's debris remediation capacity, we propose objective function: 
\begin{equation}
    \sum_{t\in {\mathcal{T}}\setminus\{T-1\}} \sum_{d\in{\mathcal{D}}}\sum_{i\in{\mathcal{J}}_{td}}\sum_{j\in{\mathcal{J}}_{t+1,d}} R_{tdij}x_{tdij} \label{eq:obj_t}
\end{equation}

We define path contiguity constraints, inspired by network flow problems, to ensure that debris $d$ occupies exactly one orbital slot at each time step $t$, and that it selects a feasible path defined in TEG ${\mathcal{H}}_d$. Constraints~\eqref{eq:flow_debris_init} require each debris to transfer from its current orbital slot $i=0$ at the epoch $t=0$, to a new one at time step $t=1$. In addition, constraints~\eqref{eq:flow_debris_eq} enforce debris to follow a contiguous path. Specifically, if at time step $t-1$ it transfers from orbital slot $i \in {\mathcal{J}}_{t-1,d}$ to orbital slot $j \in {\mathcal{J}}_{tdi}$, then at time step $t$, it has to take exactly one orbital slot $k \in {\mathcal{J}}_{t+1,dj}$. 
\begin{subequations}
    \begin{alignat}{2}
    & \sum_{j \in {\mathcal{J}}_{t+1,d}}x_{0,d,0,j} = 1,\hspace{95pt} \forall d \in \mathcal{D} \label{eq:flow_debris_init}\\
    &\sum_{k \in {\mathcal{J}}_{t+1,j}}x_{tdjk} - x_{t-1,dij} = 0, \quad \forall t \in {\mathcal{T}}\setminus\{0,T-1\}, \notag \\
    & \hspace{83pt} \forall d \in {\mathcal{D}}, \forall i \in {\mathcal{J}}_{t-1,d}, \forall j \in {\mathcal{J}}_{tdi} \label{eq:flow_debris_eq}
\end{alignat}
\end{subequations}

As illustrated in Figure~\ref{fig:debris_teg}, debris transfers to a new orbit as a consequence of at least one L2D engagement. Constraints~\eqref{eq:coupling_x_y} realize this condition by enforcing that the necessary platforms located at their corresponding orbital slots, defined in set ${\mathcal{P}}_{tdj}$, engage debris $d$ at time step $t$ such that it transfers from orbital slot $i \in {\mathcal{J}}_{td}$ to $j \in {\mathcal{J}}_{t+1,di}$.
\begin{multline}
    \sum_{(p,s) \in {\mathcal{(PS)}}_{tdj}}y^p_{tsd} \ge (PS)_{tdj} x_{tdij}, \quad \forall t \in{\mathcal{T}}\setminus\{T-1\},\\
    \quad  \forall d\in{\mathcal{D}}, \forall i \in {\mathcal{J}}_{td},\forall j \in {\mathcal{J}}_{t+1,di} \label{eq:coupling_x_y}
\end{multline}
Further, we introduce coupling constraints~\eqref{eq:coupling_y_z} to enforce that each platform performs an L2D engagement if and only if it occupies the required orbital slot.
\begin{multline}
    \sum_{d \in {\mathcal{D}}}y^p_{tsd} \le \sum_{w \in {\mathcal{S}}_{t+1,p}}z^p_{tsw}, \quad \forall t \in {\mathcal{T}}\setminus\{T-1\},\\
     \forall p\in{\mathcal{P}}, \forall s\in{\mathcal{S}}_{tp}\label{eq:coupling_y_z}
\end{multline}

We propose a set of assignment constraints to ensure the feasible operation of each platform. First, each platform can engage at most one debris per time step, as enforced by constraints~\eqref{eq:assignment_y_cardinality}. Second, each platform cannot concurrently transfer to a new orbital slot and perform an L2D engagement over debris at time step $t$; to enforce this condition, we introduce constraints~\eqref{eq:assignment_z_cardinality}.
\begin{subequations}
    \begin{alignat}{2}
    &\sum_{s \in {\mathcal{S}}_{tp}}\sum_{d \in {\mathcal{D}}} y^p_{tsd} \le 1, \quad \forall t \in {\mathcal{T}}\setminus\{T-1\}, \forall p \in {\mathcal{P}} \label{eq:assignment_y_cardinality}\\
    &\sum_{s \in {\mathcal{S}}_{tp}}\sum_{d \in {\mathcal{D}}} y^p_{tsd} + \sum_{s \in {\mathcal{S}}_{tp}}\sum_{(w\neq s) \in {{\mathcal{S}}_{t+1,p}}}z^p_{lsw} \le 1, \notag\\
    &\hspace{70pt}\forall t \in {\mathcal{T}}\setminus\{T-1\}, \forall p \in {\mathcal{P}} \label{eq:assignment_z_cardinality}
\end{alignat}
\end{subequations}

Path contiguity constraints are further defined for each platform $p$. At the epoch $t=0$, constraints~\eqref{eq:flow_platform_init} enforce each platform $p$ to transfer to a new orbital slot $w \in {\mathcal{S}}_{1,p}$. In addition, constraints~\eqref{eq:flow_platform_eq} guarantee equilibrium in the transfer path of the platform by enforcing it to select a transfer from orbital slot $s \in {\mathcal{S}}_{tp}$ to a single orbital slot $w \in {\mathcal{S}}_{t+1,p}$, only if it currently occupies orbital slot $s$.
\begin{subequations}
    \begin{alignat}{2}
    & \sum_{w \in {\mathcal{S}}_{1,p}}z^p_{0,0,w}=1,\hspace{60pt} \forall p \in {\mathcal{P}}\label{eq:flow_platform_init}\\
    &\sum_{w \in {\mathcal{S}}_{t+1,p}}z^p_{tsw} - \sum_{b \in {\mathcal{S}}_{t-1,p}}z^p_{t-1,bs} = 0, \quad  \notag \\
    & \hspace{20pt} \forall t \in {\mathcal{T}}\setminus\{0,T-1\}, \forall p \in {\mathcal{P}}, \forall s \in {\mathcal{S}}_{tp} \label{eq:flow_platform_eq}
\end{alignat}
\end{subequations}
Budget constraints are included in the optimization to impose an upper-bound on the platform's $\Delta v$ dedicated to reconfiguration.
\begin{multline}
    \sum_{t \in {\mathcal{T}}\setminus\{T-1\}}\sum_{s \in {\mathcal{S}}_{tp}}\sum_{w \in {\mathcal{S}}_{t+1,p}}c^p_{tsw}z^p_{tsw} \le c^p_{\max}, \\\forall p \in {\mathcal{P}}
\end{multline}

Lastly the binary domain for all decision variables is enforced by:
\begin{subequations}
    \begin{alignat}{2}
    &x_{tdij} \in \{0,1\}, \quad \forall t \in {\mathcal{T}}\setminus\{T-1\}, \notag\\
    &\hspace{80pt}\forall d \in {\mathcal{D}}, \forall i \in {\mathcal{J}}_{td}, \forall j \in {\mathcal{J}}_{t+1,di}\label{eq:domain_x}\\
    &y^p_{tsd} \in \{0,1\}, \quad \forall t \in {\mathcal{T}}\setminus\{T-1\}, \notag\\
    &\hspace{100pt}\forall p \in {\mathcal{P}},\forall s \in {\mathcal{S}}_{tp}, \forall d \in {\mathcal{D}}\label{eq:domain_y}\\
    &z^p_{tsw} \in \{0,1\}, \quad \forall t \in {\mathcal{T}}\setminus\{T-1\}, \notag\\
    &\hspace{80pt}\forall p \in  {\mathcal{P}},\forall s \in {\mathcal{S}}_{tp}, \forall w \in {\mathcal{S}}_{t+1,p}\label{eq:domain_z}
\end{alignat}
\end{subequations}

Piecing all constraints and the objective function, the R-L2D-ESP that maximizes the debris remediation capacity for a given mission planning horizon is given as:
\begin{alignat*}{2}
\max \quad & \text{Objective function~\eqref{eq:obj_t}} \\
\text{s.t.} \quad & \text{Constraints \eqref{eq:flow_debris_init} -- \eqref{eq:domain_z}}
\end{alignat*}

\subsection{Receding Horizon Scheduler}

The combinatorial nature of the R-L2D-ESP, owed to the binary domain of its decision variables, makes it a Nondeterministic Polynomial (NP)-hard optimization problem. In addition, as illustrated in Figs.~\ref{fig:platform_teg} and~\ref{fig:debris_teg}, the solution space exponentially expands when the number of time steps, platforms, debris, and orbital slots increase. Consequently, the problem becomes computationally intractable for large mission scenarios. 

To mitigate the rapid expansion in the solution space, we tackle the R-L2D-ESP with an RHS. In essence, the RHS is a recursive approach that schedules the actions for the current time step based on information from future time steps. Let \texttt{RHS(L,t)} be an RHS subproblem that solves the R-L2D-ESP for a receding horizon $L$, and implements the control actions for time step $t$, defined as the control horizon. These actions are the platform's transfer decision variables $\bm{z}$, and the L2D engagement scheduling decision variables $\bm{y}$. Debris transfer decision variables $\bm{x}$, which cannot be directly controlled, update the debris state vector according to the control actions scheduled $\bm{y}$ by the R-L2D-ESP.

\begin{figure}[htpb]
    \centering
    \includegraphics[width=\linewidth]{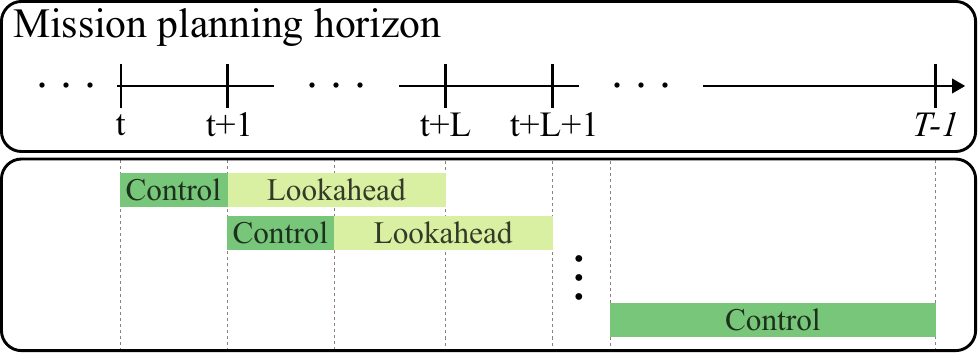}
    \caption{Illustration of the RHS approach taken to solve the R-L2D-ESP.}
    \label{fig:rh}
\end{figure}

Figure~\ref{fig:rh} illustrates the RHS used to tackle the R-L2D-ESP. Each \texttt{RHS(L,t)} is solved, control actions are implemented at the current time step $t$, and the platforms and debris states are updated for time step $t+1$. Subsequently, the receding horizon is shifted until $t=T-1-L$, and the problem is solved again. The last instance \texttt{RHS(L,T-1-L)} is solved as previously, with the exception that control actions are implemented for the remainder of the mission planning horizon.

We formally define an \text{\hypertarget{rh}{\texttt{RHS(L,t)}}} with a receding horizon length $L$ for time step $t$ as:
{\allowdisplaybreaks
\begin{subequations}
    \begin{align}
     &\max \quad \sum_{l\in \{t,\dots,t+L-1\}} \sum_{d\in{\mathcal{D}}}\sum_{i\in{\mathcal{J}}_{ld}}\sum_{j\in{\mathcal{J}}_{l+1,d}} R_{ldij}x_{ldij} \label{rhs:obj_t}\\ 
    &\sum_{j \in {\mathcal{J}}_{l+1,d}}x_{ld,0,j} = 1, \hspace{93pt}\forall d \in \mathcal{D} \label{rhs:flow_debris_init}\\
    &\sum_{k \in {\mathcal{J}}_{l+1,j}}x_{ldjk} - x_{l-1,dij} = 0, \notag \\
    &\hspace{40pt}\forall l \in \{t+1,\dots,t+L-1\}, \forall d \in {\mathcal{D}},  \notag \\
    &\hspace{116pt}\forall i \in {\mathcal{J}}_{l-1,d},  \forall j \in {\mathcal{J}}_{ldi} \label{rhs:flow_debris_eq}\\
    &\sum_{(p,s) \in {\mathcal{(PS)}}_{ldj}}y^p_{lsd} \ge (PS)_{ldj} x_{ldij},\hspace{5pt} \forall l \in \{t,\dots,t+L-1\},\notag \\
    &\hspace{80pt}\forall d\in{\mathcal{D}}, \forall i \in {\mathcal{J}}_{ld},\forall j \in {\mathcal{J}}_{l+1,di} \label{rhs:coupling_x_y}\\
    &\sum_{d \in {\mathcal{D}}}y^p_{tsd} \le \sum_{w \in {\mathcal{S}}_{t+1,p}}z^p_{tsw},\hspace{5pt} \forall l \in \{t,\dots,t+L-1\},\notag \\
    &\hspace{140pt}\forall p\in{\mathcal{P}}, \forall s\in{\mathcal{S}}_{tp}\label{rhs:coupling_y_z}\\
    &\sum_{s \in {\mathcal{S}}_{lp}}\sum_{d \in {\mathcal{D}}} y^p_{lsd} \le 1, \hspace{5pt}\forall l \in \{t,\dots,t+L-1\}, \notag \\
    &\hspace{180pt}\forall p \in {\mathcal{P}} \label{rhs:assignment_y_cardinality}\\
    &\sum_{s \in {\mathcal{S}}_{lp}}\sum_{d \in {\mathcal{D}}} y^p_{lsd} + \sum_{s \in {\mathcal{S}}_{lp}}\sum_{(w\neq s) \in {{\mathcal{S}}_{l+1,p}}}z^p_{lsw} \le 1, \notag\\
    &\hspace{70pt}\forall l \in \{t,\dots,t+L-1\}, \forall p \in {\mathcal{P}} \label{rhs:assignment_z_cardinality}\\
    & \sum_{w \in {\mathcal{S}}_{l+1,p}}z^p_{lsw}=1,\hspace{5pt}\forall p \in {\mathcal{P}}, \forall s \in \tilde{\mathcal{S}}_{lp}\label{rhs:flow_platform_init}\\
    &\sum_{w \in {\mathcal{S}}_{l+1,p}}z^p_{lsw} - \sum_{b \in {\mathcal{S}}_{l-1,p}}z^p_{l-1,bs} = 0,  \notag \\
    &\hspace{20pt}\forall l \in \{t+1,\dots,t+L-1\}, \forall p \in {\mathcal{P}}, \forall s \in {\mathcal{S}}_{lp} \label{rhs:flow_platform_eq}\\
    &\sum_{l \in \{t,\dots,t+L-1\}}\sum_{s \in {\mathcal{S}}_{lp}}\sum_{w \in {\mathcal{S}}_{l+1,p}}c^p_{lsw}z^p_{lsw} \le c^p_{\max}, \notag \\
    &\hspace{180pt}\forall p \in {\mathcal{P}} \label{rhs:platform_budget}\\
    &x_{ldij} \in \{0,1\}, \hspace{5pt} \forall l \in \{t+1,\dots,t+L-1\}, \notag\\
    &\hspace{80pt}\forall d \in {\mathcal{D}}, \forall i \in {\mathcal{J}}_{ld}, \forall j \in {\mathcal{J}}_{l+1,di}\label{rhs:domain_x}\\
    &y^p_{lsd} \in \{0,1\}, \hspace{5pt} \forall l \in \{t+1,\dots,t+L-1\}, \notag\\
    &\hspace{85pt}\forall p \in {\mathcal{P}},\forall s \in {\mathcal{S}}_{lp}, \forall d \in {\mathcal{D}}\label{rhs:domain_y}\\
    &z^p_{lsw} \in \{0,1\},\hspace{5pt} \forall l \in \{t+1,\dots,t+L-1\}, \notag\\
    &\hspace{80pt}\forall p \in  {\mathcal{P}},\forall s \in {\mathcal{S}}_{lp}, \forall w \in {\mathcal{S}}_{l+1,p}\label{rhs:domain_z}
    \end{align}
\end{subequations}}
The \hyperlink{rh}{\texttt{RHS(L,t)}} is comparable to the R-L2D-ESP formulation. Objective function~\eqref{rhs:obj_t} maximizes the debris remediation capacity for the receding horizon $L$. Constraints~\eqref{rhs:flow_debris_init} and~\eqref{rhs:flow_debris_eq} are the debris path contiguity constraints. Constraints~\eqref{rhs:coupling_x_y} are the coupling constraints between L2D engagements and debris transfer. Constraints~\eqref{rhs:coupling_y_z},~\eqref{rhs:assignment_y_cardinality}, and~\eqref{rhs:assignment_z_cardinality} are the assignment constraints. Further, constraints~\eqref{rhs:flow_platform_init}, and~\eqref{rhs:flow_platform_eq} are the platform's transfer path contiguity constraints. Since the horizon is shifted for each time step, we introduce set $\tilde{\mathcal{S}}_{lp}$ with index $s$ and cardinality $\tilde{S}_{lp}$, to define platform's $p$ unique orbital slot defined at the initial time step $l=t$. Constraints~\eqref{rhs:platform_budget} impose a $\Delta v$ transfer budget on the platform's orbital maneuverability. Lastly, constraints~\eqref{rhs:domain_x} to~\eqref{rhs:domain_z} define the decision variables' binary domain for the receding horizon.

The \hyperlink{rh}{\texttt{RHS(L,t)}} is recursively solved using Algorithm~\ref{alg:rh}. The algorithm starts by initializing the current orbital slots of platform $p$ at time step $t=0$ with $\tilde{\mathcal{S}}_{0,p}$. Then, it iterates over time steps $l=0$ to $l=T-L-2$ and performs the following sequence. First, using Algorithm~\ref{alg:debris_teg}, it generates the necessary parameters for the \hyperlink{rh}{\texttt{RHS(L,t)}}. Second, leveraging these parameters, it solves the \hyperlink{rh}{\texttt{RHS(L,t)}} for the receding horizon $\{t,\dots,t+L\}$. The optimization results are the control actions $\bm{y}(l),\bm{z}(l)$, as well as the debris transfer variables $\bm{x}(l)$, and the debris remediation capacity $V(l)$. Third, leveraging the decision variables, \texttt{UpdatePlatformBudget} is used to compute the available $\Delta v$ dedicated to platform transfers for the rest of the mission planning horizon by:
\begin{equation}
    c^p_{\max} = c^p_{\max}-\sum_{s \in {\mathcal{S}}_{lp}}\sum_{w \in {\mathcal{S}}_{l+1,p}}c^p_{lsw}z^p_{lsw}, \qquad \forall p \in {\mathcal{P}}\label{eq:update_dv_platform}
\end{equation}
Further, the set of initial platforms' positions $\tilde{\mathcal{S}}_{l+1,p}$ for the subsequent time step $l+1$ is updated with transfer decision variables $\bm{z}(l)$. Similarly, \texttt{UpdatePlatformBudget} uses debris transfer decision variables $\bm{x}(l)$ to update debris state vector for the remainder of the mission planning horizon. Lastly, defining $\bm{\alpha}\coloneqq\{T-1-L,\dots,T-1\}$, we can set the control and prediction horizons to be equal (as shown in the last block of Figure~\ref{fig:rh}), and repeat the process described above. However, the difference stems from the fact that the control actions, as well as debris and platform state updates, are implemented for the receding horizon defined by $\bm{\alpha}$.
\begin{algorithm}[htpb]
\scriptsize
    \caption{Receding Horizon Scheduler.}
    \label{alg:rh}
    \begin{algorithmic}[1]
        \renewcommand{\algorithmicrequire}{\textbf{Input:}}
        \renewcommand{\algorithmicensure}{\textbf{Output:}}
        \Require ${\mathcal{T}}, L, {\mathcal{P}}, {\mathcal{S}}^p,{\mathcal{D}}, \bm{c}_{\max}$
        \Ensure $V, \bm{x},\bm{y},\bm{z}, \bm{q}_{d}$
        \State $\tilde{\mathcal{S}}_{0,p} \gets {\mathcal{S}}_{0,p}$
        \For{$l \in \{0,\dots,T-L-2\}$}
            \State ${\mathcal{J}}_d, R_{ldij}, {\mathcal{J}}_{ldi}, {\mathcal{(PS)}}_{ldj}\gets\texttt{Algorithm}~\ref{alg:debris_teg}$ 
            \State $\bm{x}(l),\bm{y}(l),\bm{z}(l), V(l)\gets \hyperlink{rh}{\texttt{RHS(L,l)}}({\mathcal{J}}_d, R_{ldij}, {\mathcal{J}}_{ldi}, {\mathcal{(PS)}}_{ldj})$ 
            \State $c^p_{\max}\gets \texttt{UpdatePlatformBudget}(\bm{z}(l),\bm{c}_{\max})$
            \State $\tilde{\mathcal{S}}_{l+1,p}\gets\bm{z}(l)$
            \State $\bm{q}_{l+1:T-1,d} \gets \texttt{UpdateDebrisState}(\bm{x}(l))$
        \EndFor
        \State $\bm{\alpha}\gets \{T-1-L,\dots,T-1\}$
        \State ${\mathcal{J}}_d, R_{\bm{\alpha}dij}, {\mathcal{J}}_{\bm{\alpha}di}, {\mathcal{P}}_{\bm{\alpha}dj}\gets\texttt{Algorithm}~\ref{alg:debris_teg}$ 
        \State $\bm{x}(\bm{\alpha}),\bm{y}(\bm{\alpha}),\bm{z}(\bm{\alpha}), V(\bm{\alpha})\gets \hyperlink{rh}{\texttt{RHS(L,$\alpha$)}}({\mathcal{J}}_d, R_{\bm{\alpha}dij}, {\mathcal{J}}_{\bm{\alpha}di}, {\mathcal{P}}_{\bm{\alpha}dj})$
        \State $\bm{q}_{\bm{\alpha},d} \gets \texttt{UpdateDebrisState}(\bm{x}(\bm{\alpha}))$
    \end{algorithmic}
\end{algorithm}

\section{R-L2D-ESP Validation} \label{sec:case_studies_1}
We conduct a set of computational experiments to demonstrate the applicability of the R-L2D-ESP. We compare the debris remediation capacity and the number of debris deorbited for three distinct CONOPS. The first CONOPS, referred to as \textit{baseline}, assumes a non-reconfigurable constellation configuration. The second CONOPS, denoted as \textit{plane change}, consists of platforms performing plane change and phasing maneuvers to reconfigure the constellation configuration. Similarly, the third CONOPS, referred to as \textit{altitude change}, performs constellation reconfiguration by leveraging the platform's altitude change and phasing maneuvers. Alternatively, a CONOPS combining plane and altitude change with phasing could be adopted; however, it would significantly increase the number of candidate orbital slots, and consequently, the problem's complexity and solution space.

The epoch is designated as August 1, 2025, at 12:00:00.000 Coordinated Universal Time, and the mission planning horizon is set to two days, uniformly discretized in time steps of \SI{3}{min} size. Throughout the time step's length, it is assumed that each platform engages debris for a duration of at most \SI{10}{s}, and the remainder is dedicated to cooling the system. Without loss of generality, each platform utilizes the L'ADROIT space-based laser system \cite{phipps_AA_2014_ladroit}, with its key parameters summarized in Table~\ref{table:ladroit}. However, it is worth noting that the R-L2D-ESP is agnostic to the selected laser system. 
\begin{table}[htpb]
\renewcommand{\arraystretch}{1.3}
\caption{\textbf{L'ADROIT system parameters.}}
\label{table:ladroit}
\centering
\begin{tabular}{|c|c|c|}
\hline
\bfseries Parameter & \bfseries Value & \bfseries Unit \\
\hline\hline
$u_{\max}$ & 325 & \si{km}\\
$u_{\min}$ & 175 & \si{km}\\
$\phi$ & 1.5 & \si{m}\\
$M^2$ & 2 & -\\
$a$ & 1.27 & -\\
$\lambda$ & 355 &\si{nm}\\
$c_{\text{m}}$  & 100 & \si{N/MW} \\
$E$ & 380 & \si{J}\\
$\eta$ & 0.5 & -\\
$n_{\text{L2D}}$ & 560 & -\\
\hline
\end{tabular}
\end{table}

A small debris population of 395 objects with aluminum surface material composition and a surface mass density \SI{0.2}{km/m^2}, is considered in the computational experiments. Particularly, we assign each debris to a circular orbit. Debris' semi-major axis is randomly assigned using the altitude's relative frequency obtained from the European Space Agency’s Meteoroid And Space debris Terrestrial Environment Reference (MASTER-8) model \cite{MASTER8}. Conversely, the right ascension of the ascending node (RAAN) and argument of latitude are uniformly discretized between 0 and \SI{360}{deg}, and their inclination between 0 and \SI{180}{deg}.

\begin{figure}
    \centering
    \includegraphics[width=0.87\linewidth]{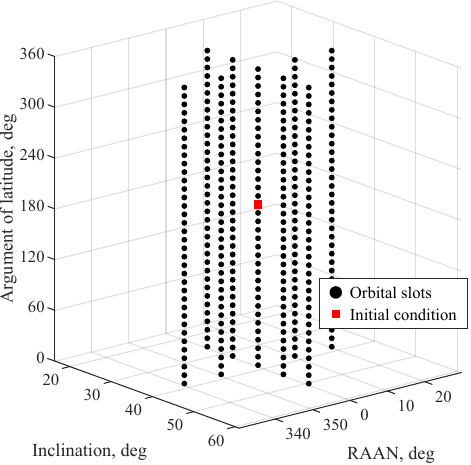}
    \caption{Plane change orbital slots for platform $p=1$ defined at the epoch.}
    \label{fig:slots_plane}
\end{figure}

\begin{figure}
    \centering
    \includegraphics[width=0.87\linewidth]{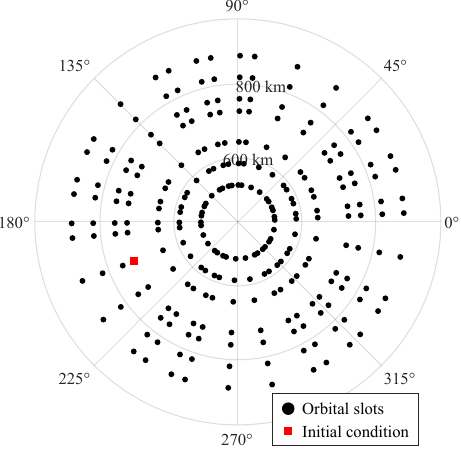}
    \caption{Altitude change orbital slots' semi-major axis versus argument of latitude plane for platform $p=1$ defined at the epoch.}
    \label{fig:slots_sma}
\end{figure}

The three CONOPS share the same initial constellation configuration consisting of six platforms. Their orbits are determined using the maximal covering location problem (MCLP) formulation proposed in Ref.~\cite{williams_asr_2025}. In essence, MCLP determines the constellation configuration that uses exactly $P$ platforms, and maximizes the L2D engagement opportunities given a user-defined debris field and mission planning horizon. Furthermore, each platform is assigned a $\Delta v$ transfer budget of \SI{2}{km/s}. The plane change CONOPS generates its candidate orbital slots by assigning an argument of latitude obtained from uniformly phasing the orbital plane in 36 steps, and generates 5 different orbital planes by changing their inclination and RAAN, as illustrated in Figure~\ref{fig:slots_plane}. These changes are determined based on the $c^p_{\max}$ budget using the rearranged boundary value problem from Ref.~\cite{vallado} with a scaling factor of 80\%. For a more detailed explanation on how to leverage the rearranged boundary value problem, we refer the reader to Appendix~\ref{app:delta_v}. Similarly, the altitude change CONOPS generates its orbital slots adopting the same phasing pattern as the one described above, and defining seven altitude layers with a step of \SI{50}{km}, where three of them increase the altitude, and three of them decrease it, as illustrated in Figure~\ref{fig:slots_sma}. 

\subsection{Results}
Table~\ref{table:opt_results} presents the debris remediation capacity for the three CONOPS and the percentage improvement over the baseline. The baseline obtains a debris remediation capacity of \num{20890.17}, and the plane change attains a debris remediation capacity of \num{27617.32}, signifying an improvement of \SI{32.20}{\%}. Similarly, the altitude change achieves a debris remediation capacity of \num{27693.67}, revealing an improvement of \SI{32.56}{\%}. In essence, the two reconfiguration-based CONOPS deliver similar improvements over the baseline, proving the value of constellation reconfiguration for orbital debris remediation.
\begin{table}[htpb]
\renewcommand{\arraystretch}{1.3}
\caption{\textbf{Optimization results.}}
\label{table:opt_results}
\centering
\begin{tabular}{|c|c|c|}
\hline
\multirow{2}{*}{\bfseries CONOPS} & \bfseries Debris remediation & \multirow{2}{*}{\bfseries Improvement}\\
 & \bfseries capacity &\\
\hline\hline
\bfseries Baseline & 20,890.17 & -\\
\bfseries Plane change &27,617.32 & 32.20 \%\\
\bfseries Altitude change  &\cellcolor[HTML]{e0ecf4}27,693.67 & \cellcolor[HTML]{e0ecf4}32.56 \%\\
\hline
\end{tabular}
\end{table}

Table~\ref{table:comparison_opt} shows the total number of debris deorbited by each CONOPS. Although useful for comparing the performance of each mission, it should be noted that this metric is distinct from the optimization objective used. The baseline constellation deorbits 104 debris objects, while the plane change and altitude change 137 and 130, respectively. The plane change and altitude change percentage improvements over the baseline case are \SI{31.73}{\%} and \SI{25}{\%}, respectively. In these numerical examples, the CONOPS that obtains the highest debris remediation capacity (\textit{i.e.,} altitude change) is distinct from the one that deorbits the most debris objects (\textit{i.e.,} plane change). This is owed to the definition of the debris remediation capacity in Eq.~\eqref{eq:reward}, which, although assigning a higher reward to the deorbiting debris, also accounts for the lowering of debris' periapsis radius closer to the deorbit altitude threshold. Therefore, throughout the mission planning horizon, the cumulative sum of lowering debris can surpass the reward collected from deorbiting debris.
\begin{table}[htpb]
\renewcommand{\arraystretch}{1.3}
\caption{\textbf{Total number of debris deorbited.}}
\label{table:comparison_opt}
\centering
\begin{tabular}{|c|c|c|}
\hline
\bfseries CONOPS& \bfseries Debris Deorbited & \bfseries Improvement\\
\hline\hline
\bfseries Baseline & 104 & - \\
\bfseries Plane change &\cellcolor[HTML]{e0ecf4}137 & \cellcolor[HTML]{e0ecf4}31.73 \%\\
\bfseries Altitude change  &130 &25 \%\\
\hline
\end{tabular}
\end{table}

\begin{figure}[htpb]
    \centering
    \includegraphics[width=0.65\linewidth]{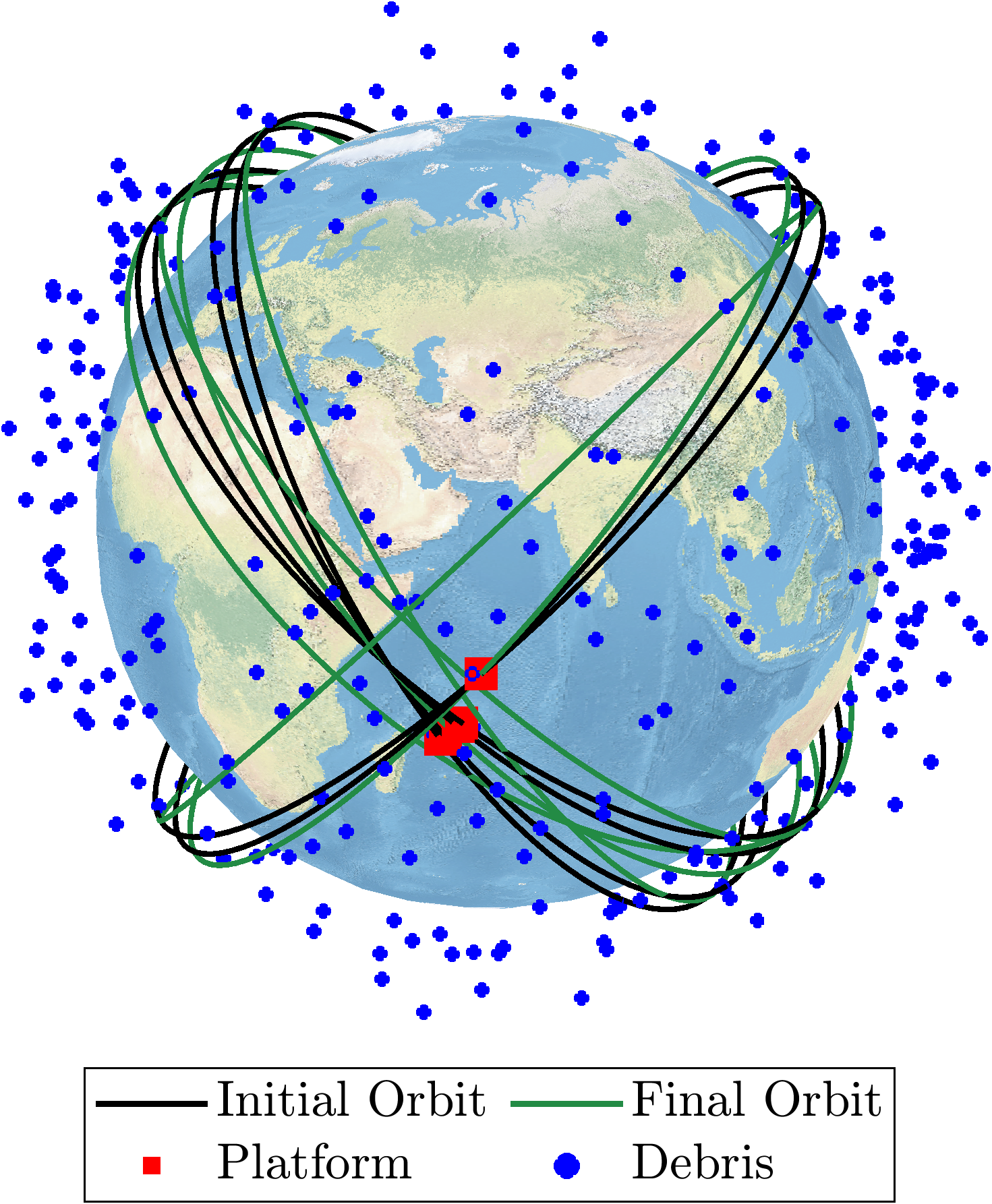}
    \caption{Plane change's initial and final constellation configurations defined at the epoch.}
    \label{fig:configuration_plane_3d}
\end{figure}

\begin{figure}[htpb]
    \centering
    \includegraphics[width=0.65\linewidth]{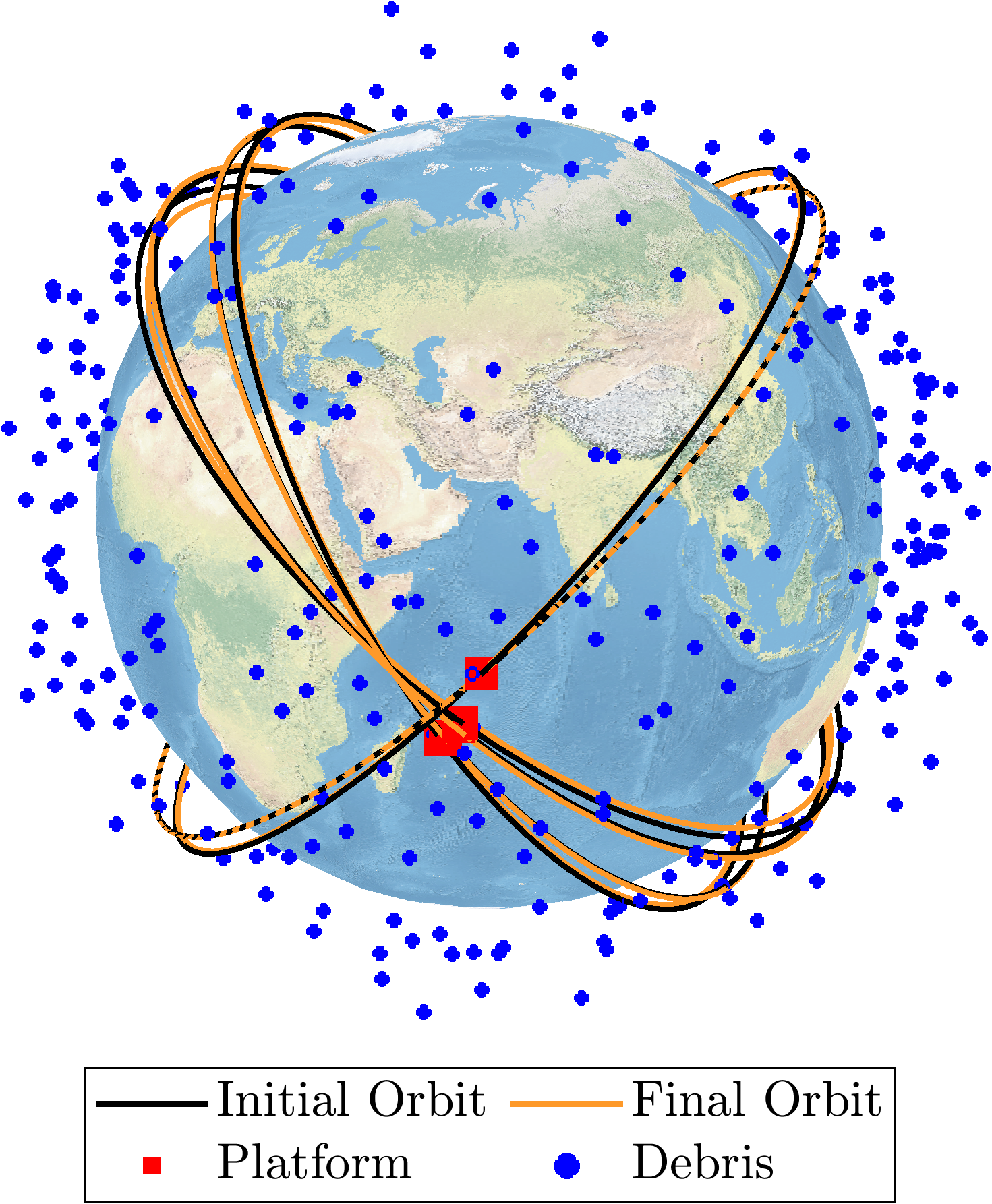}
    \caption{Altitude change's initial and final constellation configurations defined at the epoch.}
    \label{fig:configuration_sma_3d}
\end{figure}

Figures~\ref{fig:configuration_plane_3d} and~\ref{fig:configuration_sma_3d} present the plane and altitude change's initial and final constellation configurations with their orbits defined in the Earth-centered inertial (ECI) frame, respectively, and the small debris population. Throughout the mission planning horizon, the total $\Delta v$ consumed by each platform is \SI{1.99}{km}, \SI{1.99}{km}, \SI{1.98}{km}, \SI{1.99}{km}, \SI{1.99}{km}, and \SI{1.98}{km}. Similarly, the bottom configuration corresponds to the altitude change, where the final orbits are displayed in red. The total $\Delta v$ utilized is \SI{1.98}{km}, \SI{1.99}{km}, \SI{1.99}{km}, \SI{1.99}{km}, \SI{1.98}{km}, and \SI{1.98}{km}. The orbital slots corresponding to the initial and final constellation configurations are summarized in Table~\ref{table:oe_constellations}.

\begin{table*}[htpb]
\renewcommand{\arraystretch}{1.1}
\caption{\textbf{CONOPS initial and final platform's orbital elements.}}
\label{table:oe_constellations}
\centering
\begin{tabular}{|c|c|c|c|c|c|c}
\hline
&\bfseries Platform & \bfseries Semi-major axis, km & \bfseries Inclination, deg & \bfseries RAAN, deg & \bfseries Argument of latitude, deg\\
\hline\hline
\multirow{6}{*}{\makecell[c]{Initial constellation \\configuration}}         &$p=1$  & {7,104.80} & \num{38.66}& \num{0}& \num{192}\\
                &$p=2$  & {7,104.80} & \num{42.33}& \num{0}& \num{192}\\
                &$p=3$  & {7,104.80} & \num{53.33}& \num{0}& \num{192}\\
                &$p=4$  & {7,104.80} & \num{57}& \num{0}& \num{192}\\
                &$p=5$  & {7,244.80} & \num{42.33}& \num{192}& \num{0}\\
                &$p=6$  & {7,244.80} & \num{46}& \num{192}& \num{0}\\
\hline
\hline
\multirow{6}{*}{\makecell[c]{Final constellation \\configuration: \\plane change}}&$p=1$  & {7,104.80} & \num{38.66}& \num{350.12}& \num{252}\\
&$p=2$  & {7,104.80} & \num{42.33}& \num{9.14}& \num{272}\\
&$p=3$               & {7,104.80} & \num{47.19}& \num{0}& \num{2}\\
&$p=4$                       & {7,104.80} & \num{57}& \num{7.32}& \num{162}\\
&$p=5$                       & {7,244.80} & \num{48.53}& \num{192}& \num{170}\\
&$p=6$                    & {7,244.80} & \num{46}& \num{174.71}& \num{30}\\
\hline
\hline
\multirow{6}{*}{\makecell[c]{Final constellation \\configuration: \\altitude change}} &$p=1$          & {6904.80} & \num{38.66}& \num{0}& \num{52}\\
&$p=2$      & {7196.47} & \num{42.33}& \num{0}& \num{152}\\
&$p=3$                       & {7138.14} & \num{53.33}& \num{0}& \num{152}\\
&$p=4$                               & {7021.47} & \num{57}& \num{0}& \num{62}\\
&$p=5$                              & {7244.80} & \num{42.33}& \num{192}& \num{10}\\
&$p=6$                               & {7186.47} & \num{46}& \num{192}& \num{260}\\
\hline
\end{tabular}
\end{table*}

\section{Case Study: On-orbit Breakup Event Mitigation}\label{sec:breakup}
This section presents a case study that models an on-orbit breakup event and demonstrates the reconfigurable space-based laser constellation's responsiveness to time-sensitive events. For the presented case study, we assume two CONOPS with a receding horizon length $L=3$. The first one consists of a non-reconfigurable constellation configuration, referred to as \textit{baseline}, and the second CONOPS, referred to as \textit{reconfiguration}, consists of altitude change and phasing maneuvers.

The mission epoch is defined as August 1, 2025, at 12:00:00.000 Coordinated Universal Time with a mission planning horizon of one day, uniformly discretized in time steps of size \SI{180}{s}. Without loss of generality, we adopt four platforms, each equipped with an L'ADROIT laser system operating in shooting-cooling mode, as described in Section~\ref{sec:case_studies_1}.

The on-orbit breakup event is modeled as follows. Initially, we assume a single debris object on a circular orbit with a semi-major axis of \SI{6900}{km}, an inclination of \SI{60}{deg}, a RAAN of \SI{45}{deg}, and an argument of latitude of \SI{89.58}{deg}. After one hour and 28 minutes, the debris object explodes, generating 100 small debris objects of surface mass density \SI{0.2}{km/m^2}. The new debris' orbits are generated by randomly ranging the initial semi-major axis at most \SI{10}{km}, and the inclination, RAAN, and argument of latitude at most \SI{4}{deg}.

The platforms are initialized in circular orbits whose orbital parameters are presented in Appendix~\ref{app:oe_breakup}, and for reconfiguration, each one has assigned a $\Delta v$ transfer budget of \SI{1.5}{km/s}. The candidate orbital slots are generated by defining 13 altitude steps of \SI{90}{km}, which increase each platform's semi-major axis defined at the epoch, and uniformly spacing the orbital plane in 31 argument of latitude values, yielding 403 candidate orbital slots for each platform. 

\subsection{Results}

Table~\ref{table:breakup_results} presents the obtained debris remediation capacity for both CONOPS. While the baseline attains 801.26, reconfiguration obtains 9,021.26, signifying an improvement of 1,025.04~\%. Further, Table~\ref{table:breakup_comparison} outlines the number of debris deorbited, where the baseline deorbits 6 and reconfiguration 47, representing an improvement of \SI{683.33}{\%}.

\begin{table}[htpb]
\renewcommand{\arraystretch}{1.3}
\caption{\textbf{On-orbit breakup optimization results.}}
\label{table:breakup_results}
\centering
\begin{tabular}{|c|c|c|}
\hline
\multirow{2}{*}{\bfseries CONOPS} & \bfseries Debris remediation & \multirow{2}{*}{\bfseries Improvement}\\
 & \bfseries capacity &\\
\hline\hline
\bfseries Baseline & 801.86 & -\\
\bfseries Reconfiguration  &\cellcolor[HTML]{e0ecf4}9,021.26 & \cellcolor[HTML]{e0ecf4}1,025.04 \%\\
\hline
\end{tabular}
\end{table}

\begin{table}[htpb]
\renewcommand{\arraystretch}{1.3}
\caption{\textbf{Total number of debris deorbited for on-orbit breakup event.}}
\label{table:breakup_comparison}
\centering
\begin{tabular}{|c|c|c|}
\hline
\bfseries CONOPS& \bfseries Debris Deorbited & \bfseries Improvement\\
\hline\hline
\bfseries Baseline & 6 & - \\
\bfseries Reconfiguration &\cellcolor[HTML]{e0ecf4}47 & \cellcolor[HTML]{e0ecf4}683.33 \%\\
\hline
\end{tabular}
\end{table}

\begin{figure*}
    \centering
    \includegraphics[width=\linewidth]{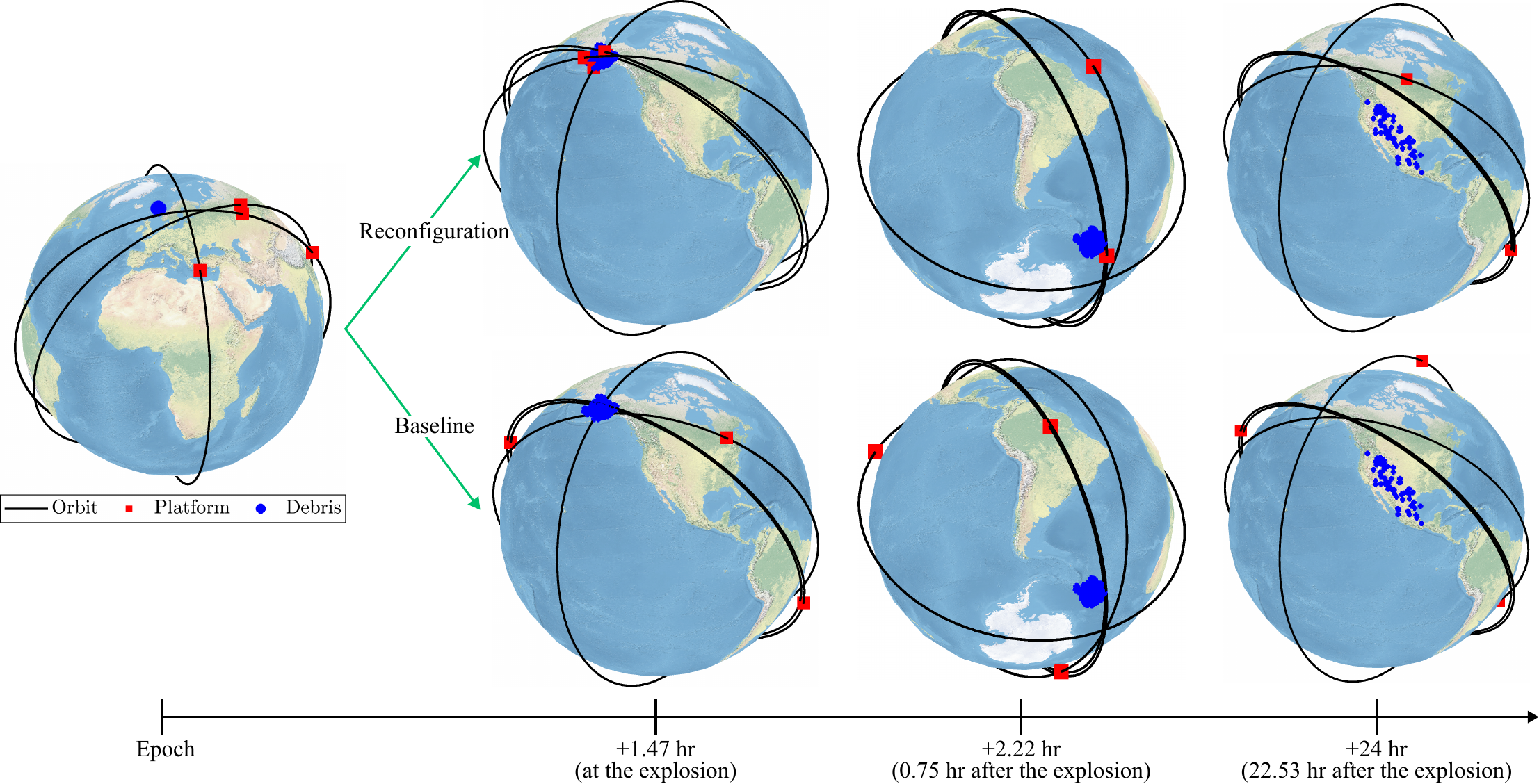}
    \caption{Timeline for baseline and reconfigurable constellations for on-orbit breakup event mitigation.}
    \label{fig:breakup_timeline}
\end{figure*}

\begin{figure}
    \centering
    \includegraphics[width=0.75\linewidth]{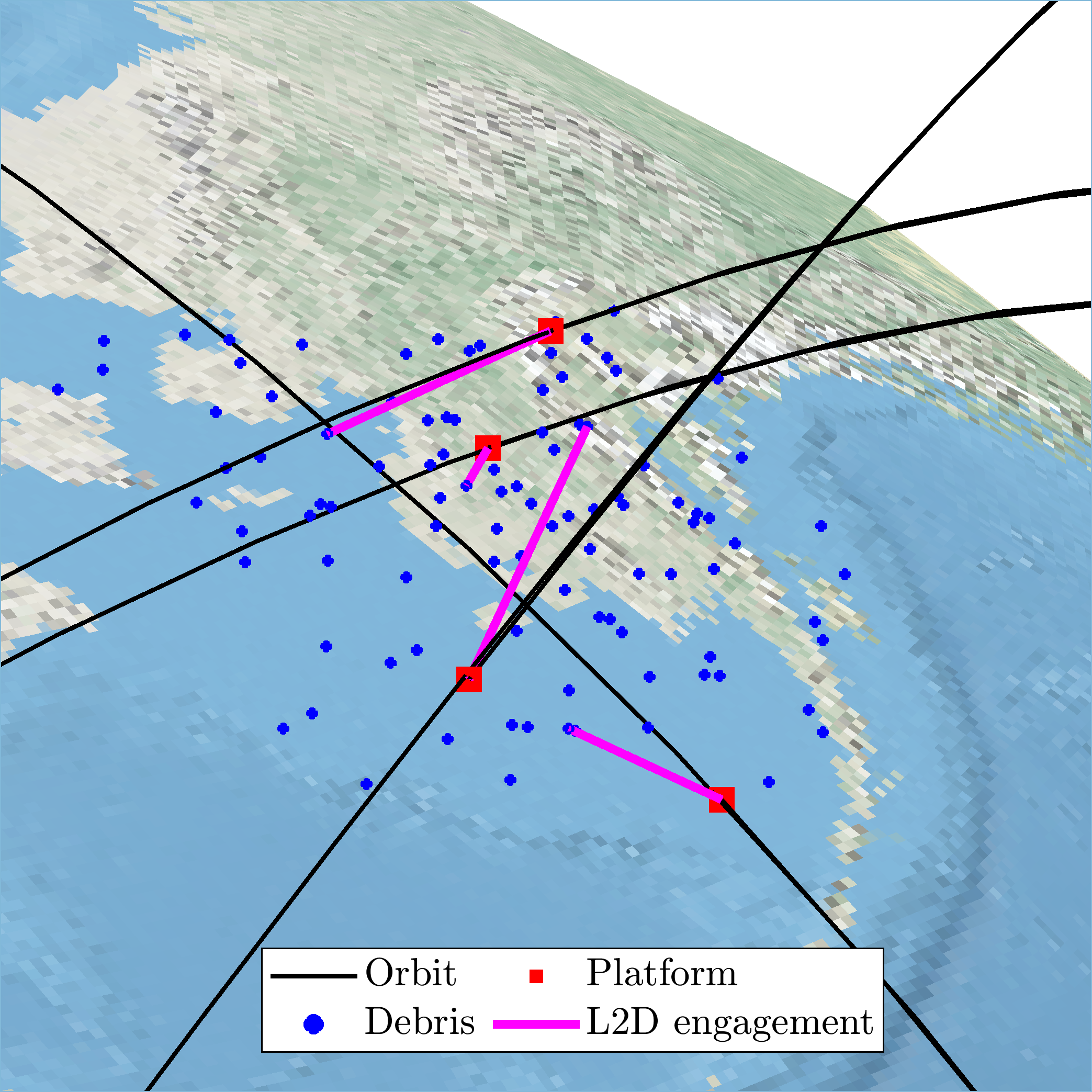}
    \caption{Reconfiguration CONOPS's constellation detailed view of L2D engagements.}
    \label{fig:breakup_zoom}
\end{figure}

Figure~\ref{fig:breakup_timeline} presents relevant events during the mission planning horizon for the proposed case study. At the epoch, both CONOPS share the same constellation configuration; however, at the time of the on-orbit breakup event, the reconfigurable CONOPS reconfigures the constellation configuration to obtain better L2D geometries. After \SI{2.22}{hr} from the epoch, the reconfigurable CONOPS is capable of tracking the debris cloud leveraging phasing maneuvers. At the end of the mission planning horizon, \textit{i.e.,} +\SI{24}{hr}, the debris cloud is spread out, clearly distinguishing the individual objects. Lastly, Figure~\ref{fig:breakup_zoom} showcases a detailed view of the reconfigurable constellation simultaneously engaging four debris objects, all of which result in instantaneous deorbit.

The case study presented in this section proves the value of constellation reconfiguration for time-sensitive events. The reconfigurable constellation significantly outperforms the non-reconfigurable one in terms of the debris remediation capacity and the number of debris deorbited. In essence, by enabling the constellation to reconfigure, a rapid and effective mitigation strategy is realized.

\section{Sensitivity Analyses}\label{sec:sensitivity}

In this section, we conduct two sensitivity analyses to characterize the impact of the receding horizon length and $\Delta v$ budget on the debris remediation capacity. Assuming the same set of parameters as in Section~\ref{sec:case_studies_1}, the first sensitivity analysis determines, for the three CONOPS, the debris remediation capacity variation for receding horizons of length $L=\{2,3,4,5\}$. The second sensitivity analysis determines the variation of the debris remediation capacity for the two reconfiguration-based CONOPS, adopting a receding horizon of length $L=3$ and establishing a platform $\Delta v$ transfer budgets of $\{0.5, 1, 1.5, 2\}$ \si{km/s} and compares them against the baseline. 

\subsection{Sensitivity Analysis I: Receding horizon length variation}

\begin{figure}[htpb]
    \centering
    \includegraphics[width=\linewidth]{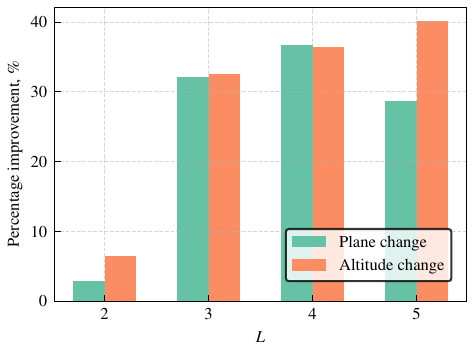}
    \caption{Debris remediation capacity percentage increase over the baseline.}
    \label{fig:percentage_increase}
\end{figure}

Table~\ref{table:l_comparison_opt} presents the debris remediation capacity for all numerical experiments, and Figure~\ref{fig:percentage_increase} shows their percentage increase with respect to the baseline. The first set of experiments corresponds to $L=2$, where the best performing CONOPS is the altitude change with a debris remediation capacity of 11,551.20. The percentage increase with respect to the baseline is \SI{6.5}{\%}. The second set of experiments corresponds to $L=3$, and has been described in Section~\ref{sec:case_studies_1}. The third set of experiments has an $L=4$ and the best-performing CONOPS is the plane change, which obtained a debris remediation capacity of 43,533.56 and presents a percentage increase of \SI{36.74}{\%} with respect to the baseline. The fourth set of numerical experiments, which adopts an $L=5$, has the altitude change as the best-performing CONOPS, with a debris remediation capacity of 59,373.73 and a percentage increase of \SI{40.09}{\%} with respect to the baseline.

\begin{table}[htpb]
\renewcommand{\arraystretch}{1.3}
\caption{\textbf{Receding horizon length versus debris remediation capacity.}}
\label{table:l_comparison_opt}
\centering
\begin{tabular}{|c|c c c|}
\hline
  \bfseries            & \multicolumn{3}{c|}{\textbf{Debris remediation capacity}}\\
\bfseries $L$ & \bfseries Baseline & \bfseries Plane change & \bfseries Altitude change\\
\hline\hline
2 & 10,846.08 & 11,167.03 & \cellcolor[HTML]{e0ecf4}11,551.20\\
3 & 20,890.17 & 27,617.32 & \cellcolor[HTML]{e0ecf4}27,693.67\\
4 & 31,836.66 & \cellcolor[HTML]{e0ecf4}43,533.56 & 43,426.79\\
5 & 42,380.88 & 54,556.66 & \cellcolor[HTML]{e0ecf4}59,373.73\\
\hline
\end{tabular}
\end{table}

\begin{figure}[htpb]
    \centering
    \includegraphics[width=\linewidth]{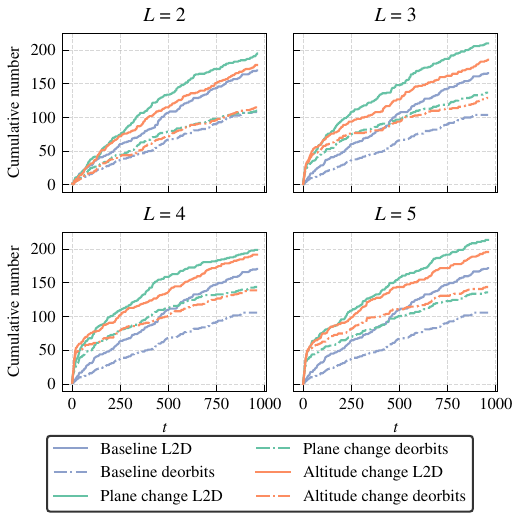}
    \caption{Cumulative number of L2D engagements and debris deorbited for different receding horizon lengths.}
    \label{fig:cumulative_comparison}
\end{figure}

Figure~\ref{fig:cumulative_comparison} presents the cumulative number of L2D engagements and debris deorbited for the four sets of numerical experiments. The plane change has the largest number of cumulative engagements for all receding horizon lengths, and for $L=3$ and $L=4$ it deorbits more debris objects than all other CONOPS, as summarized in Table~\ref{table:l_comparison_derived}. The altitude change outperforms the baseline for all receding horizon lengths and deorbits the greatest number of debris objects for $L=2$ and $L=4$.

\begin{table}[htpb]
\renewcommand{\arraystretch}{1.3}
\caption{\textbf{Receding horizon length versus number of debris deorbited}}
\label{table:l_comparison_derived}
\centering
\begin{tabular}{|c|c c c|}
\hline
     \bfseries       & \multicolumn{3}{c|}{\textbf{Total debris deorbited}}\\
\bfseries  $L$ & \bfseries Baseline & \bfseries Plane change & \bfseries Altitude change\\
\hline\hline
2 & 108 & 111 & \cellcolor[HTML]{e0ecf4}115  \\
3 & 104 & \cellcolor[HTML]{e0ecf4}137 &   130 \\
4 & 106 & \cellcolor[HTML]{e0ecf4}144 & 139\\
5 & 106 & 136 & \cellcolor[HTML]{e0ecf4}144\\
\hline
\end{tabular}
\end{table}

The results presented in this section allow us to show that increasing the receding horizon length leads to a larger debris remediation capacity and number of debris deorbited. However, it should be noted that the complexity of the problem also increases as the size of the solution space grows exponentially. Furthermore, even though the debris remediation capacity is greater for the altitude change CONOPS in three out of four instances, the absolute difference with respect to the plane change is small. Therefore, it is not possible to conclude that altitude changes will always outperform plane changes.

\subsection{Sensitivity Analysis II: $\Delta v$ transfer budget variation}
\begin{figure}[htpb]
    \centering
    \includegraphics[width=\linewidth]{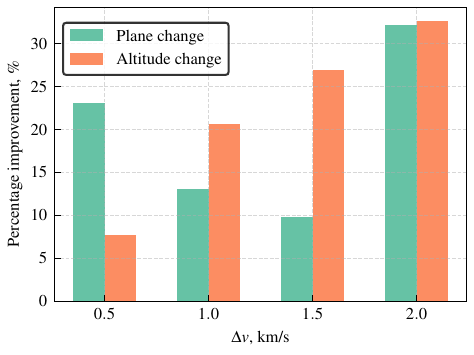}
    \caption{Debris remediation capacity percentage increase over the baseline for $L=3$.}
    \label{fig:percentage_increase_v}
\end{figure}

Table~\ref{table:v_comparison_opt} presents the debris remediation capacity for the baseline constellation and the R-L2D-ESP adopting a receding horizon length $L=3$, and different platform $\Delta v$ transfer budgets. For a transfer budget of $\Delta v=0.5$~\si{km/s}, the best performing constellation is the one with reconfigurations based on plane changes, obtaining a debris remediation capacity of 25,712.72, representing a percentage improvement over the baseline of \SI{23.08}{\%}, as illustrated in Figure~\ref{fig:percentage_increase_v}. Conversely, for transfer budgets of 1, 1.5, and \SI{2}{km/s}, the best performing CONOPS is the altitude change with debris remediation capacities of 25,207.91, 26,503.27, and 27,693.67, respectively. The percentage improvement over the baseline is \SI{20.66}{\%}, \SI{26.86}{\%}, and \SI{32.56}{\%}, respectively. 

\begin{table}[htpb]
\renewcommand{\arraystretch}{1.3}
\caption{\textbf{Per-platform $\Delta v$ budget versus debris remediation capacity.}}
\label{table:v_comparison_opt}
\centering
\begin{tabular}{|c|c c c|}
\hline
  \bfseries $\Delta v$           & \multicolumn{3}{c|}{\textbf{Debris remediation capacity}}\\
\bfseries budget & \bfseries Baseline & \bfseries Plane change & \bfseries Altitude change\\
\hline\hline
0.5 & 20,890.17 & \cellcolor[HTML]{e0ecf4}25,712.72 &  22,499.60\\
1 & 20,890.17 & 23,624.25 & \cellcolor[HTML]{e0ecf4}25,207.91\\
1.5 & 20,890.17 &  22,936.36 & \cellcolor[HTML]{e0ecf4}26,503.27\\
2 & 20,890.17 & 27,617.32& \cellcolor[HTML]{e0ecf4}27,693.67\\
\hline
\end{tabular}
\end{table}
As presented in Table~\ref{table:v_comparison_opt} and Figure~\ref{fig:percentage_increase_v}, altitude change presents a monotonic increase in its debris remediation capacity with respect to the $\Delta v$ transfer budget. Conversely, plane change obtains a higher debris remediation capacity for a $\Delta v$ transfer budget of \SI{0.5}{km/s} than for \SI{1}{km/s} and \SI{1.5}{km/s}. This stems from the receding heuristic approach adopted to tackle the R-L2D-ESP. In essence, the \hyperlink{rh}{\texttt{RHS(L,t)}} obtains a global optimum for a subproblem of length $L$; however, this solution does not consider the full length of the mission planning horizon, therefore, the scheduled constellation reconfigurations and L2D engagements, for the current subproblem, might not be the best ones for the remainder of the mission given the lack of future information. As a consequence, a larger $\Delta v$ transfer budget does not always guarantee a higher overall debris remediation capacity.

\begin{figure}[htpb]
    \centering
    \includegraphics[width=\linewidth]{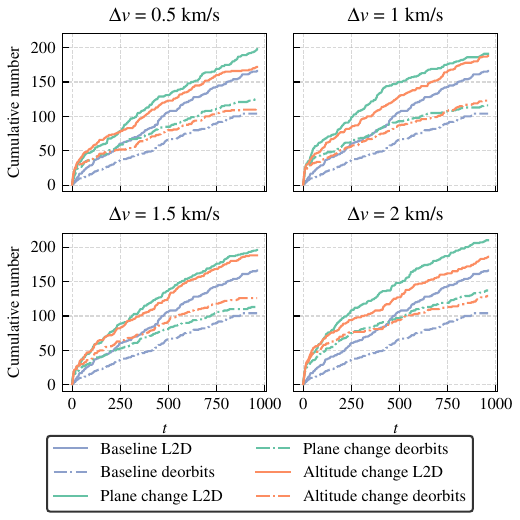}
    \caption{Cumulative number of L2D engagements and debris deorbited for different $\Delta v$ transfer budgets.}
    \label{fig:cumulative_comparison_v}
\end{figure}

Table~\ref{table:v_comparison_derived} outlines the total number of debris deorbited by each constellation. As presented in the previous section, the baseline deorbits 104 debris, while for a transfer budget of \SI{2}{km/s} the plane and altitude change deorbit 137 and 130 debris objects, respectively. Similarly, for a transfer budget of \SI{0.5}{km/s}, the plane change deorbits the largest number of debris objects, being 127. Conversely, for transfer budgets of 1 and \SI{1.5}{km/s}, the altitude change outperforms the plane change, deorbiting 124 and 126 debris objects, respectively. Figure~\ref{fig:cumulative_comparison_v} displays, in addition to the cumulative number of debris deorbited, the cumulative number of L2D engagements for each instance. All four instances exhibit the same trend, with the plane change having the largest cumulative number of L2D engagements.
\begin{table}[htpb]
\renewcommand{\arraystretch}{1.3}
\caption{\textbf{Per-platform $\Delta v$ budget versus number of debris deorbited}}
\label{table:v_comparison_derived}
\centering
\begin{tabular}{|c|c c c|}
\hline
     \bfseries $\Delta v$      & \multicolumn{3}{c|}{\textbf{Total debris deorbited}}\\
\bfseries  budget & \bfseries Baseline & \bfseries Plane change & \bfseries Altitude change\\
\hline\hline
0.5 & 104 & \cellcolor[HTML]{e0ecf4}127 &  111\\
1 & 104 & 116 & \cellcolor[HTML]{e0ecf4}124\\
1.5 & 104 &  113 & \cellcolor[HTML]{e0ecf4}126\\
2 & 104 & \cellcolor[HTML]{e0ecf4}137&  130\\
\hline
\end{tabular}
\end{table}

The numerical experiments presented in this section provide insight into the impact of the $\Delta v$ transfer budget on the debris remediation capacity. Overall, all instances significantly outperform the baseline case. For the altitude change CONOPS, there is a clear correlation between increasing the transfer budget and achieving a higher debris remediation capacity, as well as a larger number of debris deorbited. In contrast, for a plane change CONOPS, there is no correlation between increasing the transfer budget and obtaining a larger debris remediation capacity as a consequence of the receding horizon approach used to tackle the R-L2D-ESP.

\section{Conclusions}\label{sec:conclusions}
This paper proposes R-L2D-ESP, an ILP optimization problem that maximizes the debris remediation capacity by scheduling the L2D engagements and the sequence of constellation reconfigurations. We leverage the principle of PLA to model L2D engagements, and TEGs to model the orbital transfers for both debris and platforms. Given the combinatorial nature of the problem, we propose a receding horizon scheduler to solve large problem instances without suffering an explosion of the solution space.

We conduct two case studies to demonstrate the applicability of the novel R-L2D-ESP. The first case study, intended to validate the optimization problem, considers three CONOPS; the first CONOPS represents a non-reconfigurable constellation configuration. The second CONOPS performs constellation reconfiguration by leveraging phasing and plane change maneuvers, and the third one utilizes phasing and altitude change maneuvers. The two reconfigurable constellations significantly outperform the baseline for small debris remediation. The second case study considers a non-reconfigurable and a reconfigurable constellation configuration to address an on-orbit breakup event, showcasing the benefits of implementing reconfiguration against time-sensitive events. Furthermore, we conduct two sensitivity analyses to characterize the impact of the receding horizon length and $\Delta v$ budget on the debris remediation capacity. The results obtained showcase that, overall, by augmenting the receding horizon length and the $\Delta v$ transfer budget, the debris remediation capacity increases.

The case studies presented in this article prove the value of constellation reconfiguration for orbital debris remediation. First, by adopting the same number of platforms, a larger population of debris can be remediated, facilitating the scalability against the large debris population. Second, time-sensitive events can be rapidly mitigated by adopting better configurations. In essence, enabling a space-based laser constellation to reconfigure during a debris remediation mission leads to flexible and responsive mission planning. Therefore, it offers an effective and promising solution to the orbital debris problem. 

The contributions of this paper open fruitful future research directions. First, we aim to characterize the performance of different laser platform systems under the same mission scenario. Second, consider the uncertainty of the post-L2D debris trajectory due to unknown debris shape and material.

\appendices{}
\section{R-L2D-ESP Nomenclature}\label{app:nomenclature}
Table~\ref{table:r_l2d_esp} presents the symbols and their respective descriptions for the parameters and decision variables utilized to construct TEGs and in the R-L2D-ESP.
\begin{table}[htpb]
\renewcommand{\arraystretch}{1.3}
\caption{\textbf{R-L2D-ESP's parameters, and decision variables.}}
\label{table:r_l2d_esp}
\centering
\begin{tabular}{|c|l|}
\hline
\bfseries Symbol & \bfseries Description \\
\hline\hline
\multicolumn{2}{|c|}{\bfseries Parameters} \\
\hline
${\mathcal{G}}^p$  & Platform $p$'s TEG \\
${\mathcal{S}}^p$ & Set of orbital slots for platform $p$\\
${\mathcal{S}}^p_{t}$ & Set of orbital slots for platform $p$ at time\\
                    &step $t$ (indices $s,w$, cardinality $S^p_{t}$)\\
$c^p_{tsw}$ & Platform $p$'s reconfiguration cost at time\\
                    &step $t$ from orbital slot $s$ to $w$\\ 
${\mathcal{C}}_t$& Constellation configuration at time step $t$\\
${\mathcal{H}}^d$  & Debris $d$'s TEG \\
${\mathcal{J}}_d$ & Set of orbital slots for debris $d$\\
${\mathcal{J}}_{td}$ & Set of orbital slots for debris $d$ at time\\
                &step $t$ (indices $i,j$, cardinality $J_{td}$)\\
${\mathcal{(PS)}}_{t-1,dj}$ & Set of platforms $p$ at orbital slots $s$ that\\
                &engage debris $d$ at time step $t$ and generate \\
                &orbital slot $j$\\
${\mathcal{J}}_{tdi}$ & Set of orbital slots for debris $d$ at time\\
                &step $t$ that share parent orbital slot $i$ \\
                &(index $j$, cardinality $J_{tdi}$)\\
$R_{tdij}$ & Debris $d$ reward for transfering at time \\
                &step $t$ from orbital slot $i$ to $j$\\
$\Gamma_{tdij}$ & Reward $R_{tdij}$'s collision-aware penalty\\
$r_{\text{deorbit}}$ & Debris orbital radius deorbit threshold\\
$\gamma_{\text{[peri]},tdij}$ & Reward $R_{tdij}$'s periapsis radius\\
                &reduction reward\\
\hline
\hline
\multicolumn{2}{|c|}{\bfseries Decision variables} \\
\hline
$z^p_{tsw}$& Platform $p$'s dec. var. for transferring at \\
                &time step $t$ from orbital slot $s$ to $w$\\
$y^p_{tsd}$ & Platform $p$ at orbital slot $s$'s L2D \\
                &engagement dec. var. over debris $d$\\
                &at time step $t$\\
$x_{tdij}$ & Debris $d$'s dec. var. for transferring at time \\
                &step $t$ from orbital slot $i$ to $j$\\
\hline
\end{tabular}
\end{table}

\section{Rearranged Boundary Value Problem}\label{app:delta_v}
To generate the orbital slots corresponding to the plane change CONOPS, we adopt from Ref.~\cite{vallado}'s Chapter 6 the orbital boundary value problems. In essence, these problems enable us to compute the maximum change in inclination or RAAN, assuming that the platform uses the entire $\Delta v$ transfer budget scaled by a constant $\beta$. To compute the maximum change in inclination, we use:
\begin{equation}
    \delta=\Delta inc = 2\beta \text{arcsin} \left( \frac{\Delta v}{2\sqrt{\mu_\text{Earth}/r_p}} \right)
\end{equation}
where $\Delta inc$ is the maximum change in inclination, $\mu_\text{Earth}$ is the Earth's gravitational constant and $r_p$ the orbital radius of platorm $p$. Similarly, to compute the maximum change in RAAN, we use:
\begin{equation}
    \Delta \Omega = \beta \text{arccos} \left( \frac{\delta-\cos^2{(inc_p)}}{\sin^2{(inc_p)}} \right)
\end{equation}
where $inc_p$ is the inclination of platform $p$.

\section{Case Study Platform's Orbital Elements}\label{app:oe_breakup} 
Table~\ref{table:oe_breakup} summarizes the initial orbital elements for both CONOPS.
\begin{table}[htpb]
    \renewcommand{\arraystretch}{1.1}
    \caption{\textbf{CONOPS initial platform's orbital elements.}}
    \label{table:oe_breakup}
    \centering
    \begin{tabular}{|c|c|}
    \hline
    \bfseries Platform & \bfseries Semi-major axis, km\\
    \hline\hline
    $p=1$  & {6750} \\
    $p=2$  & {6750}\\
    $p=3$  & {6750}\\
    $p=4$  & {6750}\\
    \hline
    \hline
    \bfseries Platform  & \bfseries Inclination, deg\\
    \hline\hline
    $p=1$  &  \num{60}\\
    $p=2$  &  \num{87.85}\\
    $p=3$  &  \num{60}\\
    $p=4$  &  \num{60}\\
    \hline
    \hline
    \bfseries Platform & \bfseries RAAN, deg \\
    \hline\hline
    $p=1$  & \num{67.50}\\
    $p=2$  &\num{337.50}\\
    $p=3$  &\num{67.50}\\
    $p=4$  &\num{90}\\
    \hline
    \hline
    \bfseries Platform & \bfseries Argument of latitude, deg\\
    \hline\hline
    $p=1$  &\num{145.50}\\
    $p=2$  & \num{145.50}\\
    $p=3$  &  \num{109.50}\\
    $p=4$  & \num{99.50}\\
    \hline
    \end{tabular}   
\end{table}

\clearpage
\acknowledgements 
This work was supported by an Early Career Faculty grant from NASA’s Space Technology Research Grants Program under award No. 80NSSC23K1499.

\begingroup
\sloppy
\hyphenpenalty=10000
\exhyphenpenalty=10000
\bibliographystyle{IEEEtran}
\bibliography{references}

\begin{thebibliography}{10}
\providecommand{\url}[1]{#1}
\csname url@samestyle\endcsname
\providecommand{\newblock}{\relax}
\providecommand{\bibinfo}[2]{#2}
\providecommand{\BIBentrySTDinterwordspacing}{\spaceskip=0pt\relax}
\providecommand{\BIBentryALTinterwordstretchfactor}{4}
\providecommand{\BIBentryALTinterwordspacing}{\spaceskip=\fontdimen2\font plus
\BIBentryALTinterwordstretchfactor\fontdimen3\font minus \fontdimen4\font\relax}
\providecommand{\BIBforeignlanguage}[2]{{%
\expandafter\ifx\csname l@#1\endcsname\relax
\typeout{** WARNING: IEEEtran.bst: No hyphenation pattern has been}%
\typeout{** loaded for the language `#1'. Using the pattern for}%
\typeout{** the default language instead.}%
\else
\language=\csname l@#1\endcsname
\fi
#2}}
\providecommand{\BIBdecl}{\relax}
\BIBdecl

\bibitem{danzmann_asr_2003}
L.~K.~Danzmann, ``{LISA} — an {ESA} cornerstone mission for the detection and observation of gravitational waves,'' \emph{Advances in Space Research}, vol.~32, no.~7, pp. 1233--1242, 2003, available: \url{https:/doi.org/10.1016/S0273-1177(03)90323-1}.

\bibitem{nasa_iss}
{NASA}, ``{International Space Station},'' \url{https://www.nasa.gov/international-space-station/}, last accessed 26 August 2025.

\bibitem{whealan_SP_2019}
K.~{Whealan George}, ``The economic impacts of the commercial space industry,'' \emph{Space Policy}, vol.~47, pp. 181--186, 2019, available: \url{https:/doi.org/10.1016/j.spacepol.2018.12.003}.

\bibitem{faa_launch}
R.~Schaufele, M.~Lukacs, and {et al.}, ``{FAA} aerospace forecast fiscal years 2025-2045,'' FAA, Tech. Rep., 2025, \url{https://www.faa.gov/data_research/aviation/aerospace_forecasts/2025-commercial-space.pdf}.

\bibitem{colvin_nasa_2024}
J.~Locke, T.~J. Colvin, L.~Ratliff, A.~Abdul-Hamid, and C.~Samples, ``Cost and benefit analysis of mitigating, tracking, and remediating orbital debris,'' NASA, Tech. Rep., 2024, \url{https://ntrs.nasa.gov/api/citations/20240003484/downloads/2024%20-%20OTPS%20-%20CBA%20of%20Orbital%20Debris%20Phase%202%20v3.pdf}.

\bibitem{mcknight_AA_2021}
D.~McKnight, R.~Witner, F.~Letizia, S.~Lemmens, L.~Anselmo, C.~Pardini, A.~Rossi, C.~Kunstadter, S.~Kawamoto, V.~Aslanov, J.-C. {Dolado Perez}, V.~Ruch, H.~Lewis, M.~Nicolls, L.~Jing, S.~Dan, W.~Dongfang, A.~Baranov, and D.~Grishko, ``Identifying the 50 statistically-most-concerning derelict objects in {LEO},'' \emph{Acta Astronautica}, vol. 181, pp. 282--291, 2021, available: \url{https:/doi.org/10.1016/j.actaastro.2021.01.021}.

\bibitem{colvin_nasa_2023}
T.~J. Colvin, J.~Karcz, and G.~Wusk, ``Cost and benefit analysis of orbital debris remediation,'' NASA, Tech. Rep., 2023, \url{https://ntrs.nasa.gov/citations/20230002817}.

\bibitem{fang_aa_2019}
Y.~Fang, J.~Pan, Y.~Luo, and C.~Li, ``Effects of deorbit evolution on space-based pulse laser irradiating centimeter-scale space debris in {LEO},'' \emph{Acta Astronautica}, vol. 165, pp. 184--190, 2019, available: \url{https://doi.org/10.1016/j.actaastro.2019.09.010}.

\bibitem{pieters_asr_2023}
L.~Pieters and R.~Noomen, ``Simulating arbitrary interactions between small-scale space debris and a space-based pulsed laser system,'' \emph{Advances in Space Research}, vol.~72, no.~7, pp. 2778--2785, 2023, available: \url{https://doi.org/10.1016/j.asr.2022.04.049}.

\bibitem{phipps_AA_2014_ladroit}
C.~R. Phipps, ``{L'ADROIT}--a spaceborne ultraviolet laser system for space debris clearing,'' \emph{Acta Astronautica}, vol. 104, no.~1, pp. 243--255, 2014, available: \url{https://doi.org/10.1016/j.actaastro.2014.08.007}.

\bibitem{phipps_aa_2016}
C.~R. Phipps and C.~Bonnal, ``A spaceborne, pulsed uv laser system for re-entering or nudging {LEO} debris, and re-orbiting {GEO} debris,'' \emph{Acta Astronautica}, vol. 118, pp. 224--236, 2016, available: \url{https://doi.org/10.1016/j.actaastro.2015.10.005}.

\bibitem{williams_asr_2025}
D.~O. {Williams Rogers}, M.~C. Fox, P.~R. Stysley, and H.~Lee, ``Optimal placement and coordinated scheduling of distributed space-based lasers for orbital debris remediation,'' \emph{Advances in Space Research}, vol.~76, no.~9, pp. 5265--5293, 2025, available: \url{https://doi.org/10.1016/j.asr.2025.07.093}.

\bibitem{baker_asc_2025}
G.~M. Baker and H.~Lee, ``Reinforcement learning-based task planning of space-based lasers for orbital debris remediation,'' in \emph{AAS/AIAA Astrodynamics Specialist Conference}, 2025.

\bibitem{pearl_jsr_2025}
B.~D. Pearl, L.~P. Gold, and H.~Lee, ``Benchmarking agility and reconfigurability in satellite systems for tropical cyclone monitoring,'' \emph{Journal of Spacecraft and Rockets}, vol.~62, no.~4, pp. 1138--1151, 2025, available: \url{https://doi.org/10.2514/1.A36177}.

\bibitem{lee_jsr_2023}
H.~Lee and K.~Ho, ``Regional constellation reconfiguration problem: Integer linear programming formulation and lagrangian heuristic method,'' \emph{Journal of Spacecraft and Rockets}, vol.~60, no.~6, pp. 1828--1845, 2023, available: \url{https://doi.org/10.2514/1.A35685}.

\bibitem{lee_jsr_2024}
H.~Lee, D.~O. {Williams Rogers}, B.~D. Pearl, H.~Chen, and K.~Ho, ``Deterministic multistage constellation reconfiguration using integer programming and sequential decision-making methods,'' \emph{Journal of Spacecraft and Rockets}, vol.~62, no.~1, pp. 206--222, 2025, available: \url{https://doi.org/10.2514/1.A35990}.

\bibitem{pearl_arxiv_2025}
B.~D. Pearl, J.~M. Miller, and H.~Lee, ``The reconfigurable {Earth} observation satellite scheduling problem,'' \emph{Journal of Aerospace Information Systems}, 2025, (Accepted) available: \url{https://doi.org/10.48550/arXiv.2507.10394}.

\bibitem{vonderlinde_apss_2000}
D.~{von der Linde} and K.~Sokolowski-Tinten, ``The physical mechanisms of short-pulse laser ablation,'' \emph{Applied Surface Science}, vol. 154-155, pp. 1--10, 2000, available: \url{https://doi.org/10.1016/S0169-4332(99)00440-7}.

\bibitem{Stafe_LA_book_2014}
\BIBentryALTinterwordspacing
N.~N.~P. Mihai~Stafe, Aurelian~Marcu, \emph{Pulsed Laser Ablation of Solids}.\hskip 1em plus 0.5em minus 0.4em\relax Springer, 2014. [Online]. Available: \url{https://link.springer.com/book/10.1007/978-3-642-40978-3}
\BIBentrySTDinterwordspacing

\bibitem{lunney_ass_1998}
J.~G. Lunney and R.~Jordan, ``Pulsed laser ablation of metals,'' \emph{Applied Surface Science}, vol. 127-129, pp. 941--946, 1998, available: \url{https://doi.org/10.1016/S0169-4332(97)00770-8}.

\bibitem{yuan_jap_2012}
H.~Yuan, H.~Tong, M.~Li, and C.~Sun, ``Computational study of nanosecond pulsed laser ablation and the application to momentum coupling,'' \emph{Journal of Applied Physics}, vol. 112, no.~2, p. 023105, 2012, available: \url{https://doi.org/10.1063/1.4737188}.

\bibitem{soulard_AA_2014}
R.~Soulard, M.~N. Quinn, T.~Tajima, and G.~Mourou, ``{ICAN}: A novel laser architecture for space debris removal,'' \emph{Acta Astronautica}, vol. 105, no.~1, pp. 192--200, 2014, available: \url{https://doi.org/10.1016/j.actaastro.2014.09.004}.

\bibitem{MASTER8}
{European Space Agency}, ``Maintenance of the esa master model,'' \url{https://sdup.esoc.esa.int/master/downloads/documentation/7.02/MASTER_Final_Report.pdf}, June 2011.

\bibitem{vallado}
D.~Vallado, \emph{Fundamentals of Astrodynamics and Applications}.\hskip 1em plus 0.5em minus 0.4em\relax Space technology library, Microcosm Press, 2022, ch.~6, pp. 321--428.

\end{thebibliography}
\endgroup

\thebiography
\begin{biographywithpic}
{David Williams Rogers}{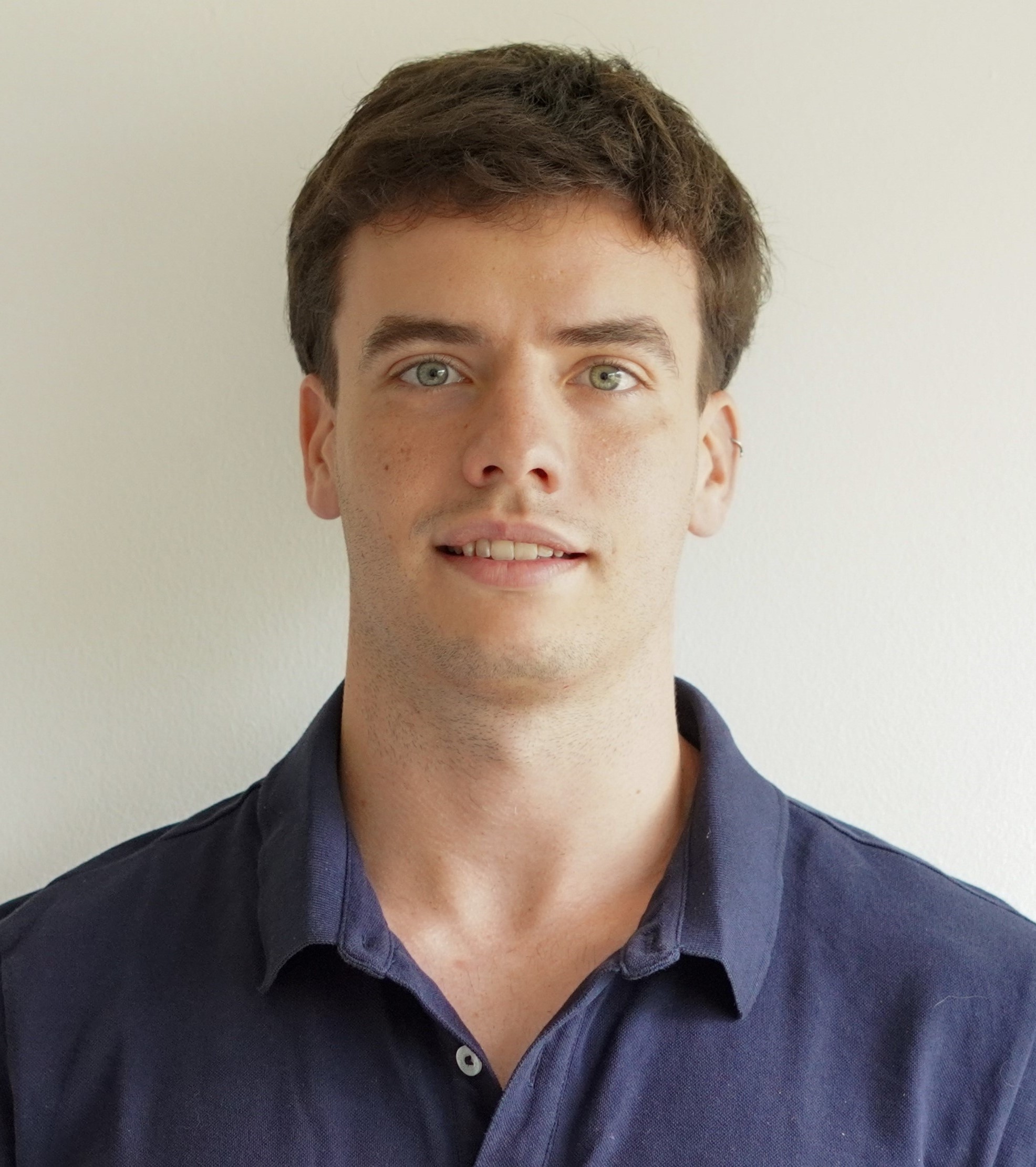}
is an aerospace engineering Ph.D. student at West Virginia University. He received his bachelor's degree in aerospace engineering from Universidad Nacional de La Plata, being the first aerospace engineer to graduate in Argentina. David is a graduate research assistant at the Space Systems Operations Research Laboratory, where his research interests lie in the intersection between distributed space systems and operations research.
\end{biographywithpic} 

\begin{biographywithpic}
{Hang Woon Lee}{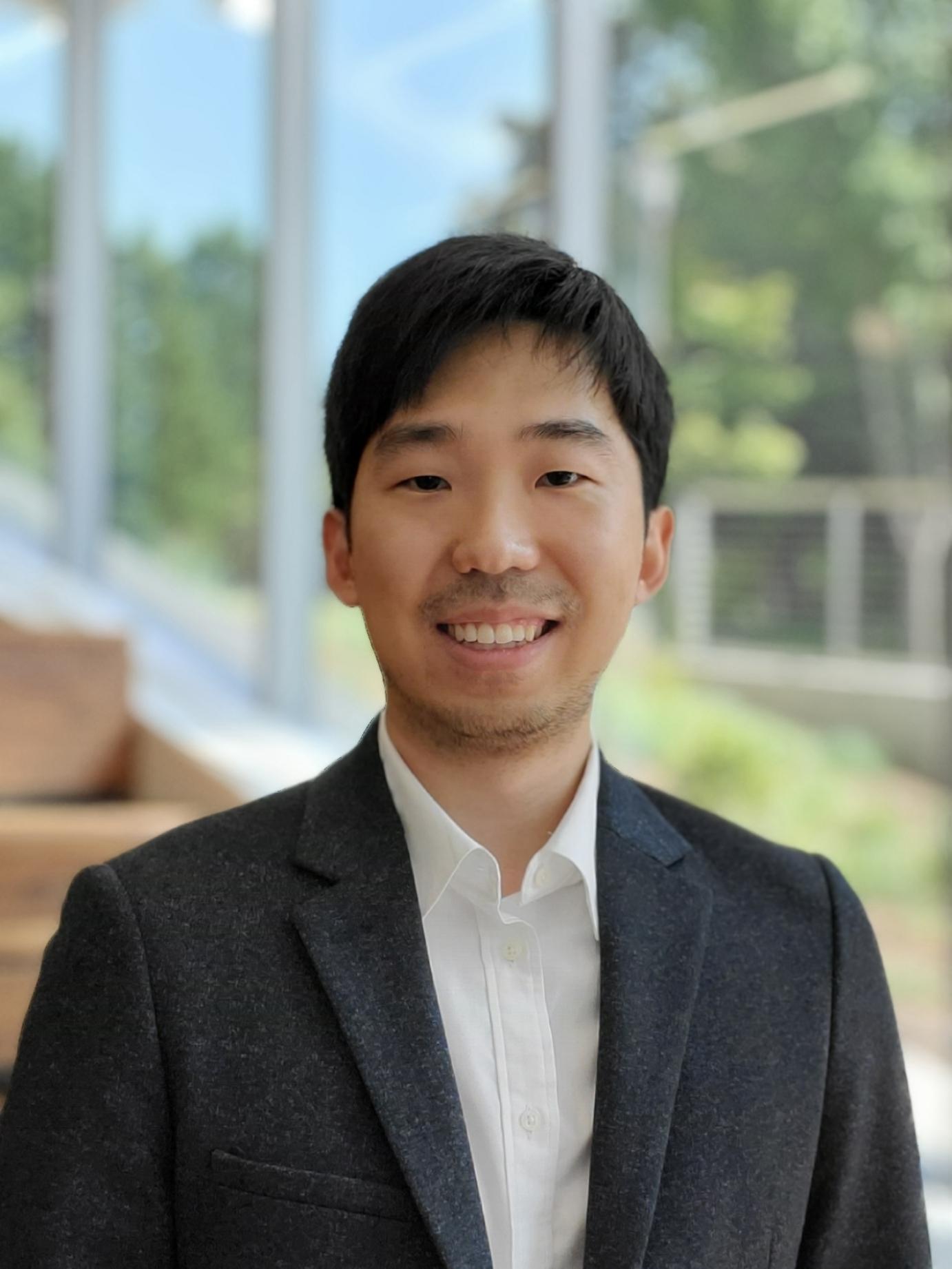}
is an assistant professor in the Department of Mechanical, Materials and Aerospace Engineering and the Herbert P. Dripps Faculty Fellow at West Virginia University. Dr. Lee directs the Space Systems Operations Research Laboratory, researching at the intersection of space systems engineering, operations research, and astrodynamics. Specifically, his research lab conducts research in the domains of distributed satellite systems, satellite networks, remote sensing, space traffic management, and space domain awareness.
\end{biographywithpic}

\end{document}